\documentclass[12pt,leqno]{amsproc}
\usepackage{amsmath,amsthm,amssymb}

\usepackage[hidelinks]{hyperref}

\usepackage{slashed}

\usepackage{xcolor}
\usepackage{url}

\usepackage{dsfont,esint}

\allowdisplaybreaks

\usepackage[all]{xy}
\xyoption{arc}

\usepackage{amsfonts}
\usepackage{enumerate}

\numberwithin{equation}{section}

 \swapnumbers

\theoremstyle{plain}
\newtheorem*{theorem*}{Theorem} 
\newtheorem*{lemma*}{Lemma}
\newtheorem*{assumption*}{Assumption}

\newtheorem{theorem}[equation]{Theorem} 
\newtheorem{lemma}[equation]{Lemma}
\newtheorem{corollary}[equation]{Corollary}
\newtheorem{proposition}[equation]{Proposition}
\newtheorem{proposition/definition}[equation]{Proposition/Definition}

\theoremstyle{definition}
\newtheorem{definition}[equation]{Definition}

\newtheorem{remark}[equation]{Remark}
\newtheorem{remarks}[equation]{Remarks}
\newtheorem*{remark*}{Remark}
\newtheorem*{remarks*}{Remarks}
\newtheorem*{observation*}{Observation}

\theoremstyle{remark}

\newcommand{\R}{\mathbb{R}}
\newcommand{\C}{\mathbb{C}}
\newcommand{\Z}{\mathbb{Z}}

\newcommand{\Tcirc}{\dot T^*}

\DeclareMathOperator{\ad}{ad}
\DeclareMathOperator{\order}{order}
\DeclareMathOperator{\SymbolAlg}{\mathrm{A}}

\DeclareMathOperator{\Tr}{Tr}
\DeclareMathOperator{\STr}{STr}
\DeclareMathOperator{\str}{str}
\DeclareMathOperator{\Todd}{Todd}

\DeclareMathOperator{\Symbol}{\operatorname{Symb}}

\DeclareMathOperator{\End}{End}

\DeclareMathOperator{\Index}{Index}

\newcommand{\bimodorder}{\operatorname{bimodule-order}}
\newcommand{\perrotorder}{\operatorname{Perrot-order}} 
\newcommand{\degree}{\operatorname{degree}}

\newcommand{\Dirac}{\slashed{D}}

\newcommand{\PDO}{\operatorname{PDO}}
\newcommand{\FPDO}{\operatorname{FPDO}}
\newcommand{\FB}{\operatorname{FB}}

\newcommand{\Poly}{\operatorname{P}}
\newcommand{\Density}{\operatorname{Dens}}
\newcommand{\gr}{\operatorname{gr}}

\newcommand{\FA}{\mathrm{FA}}
\begin{document}

\title{On Perrot's Index Cocycles}

\author{Jonathan Block}
\address{Department of Mathematics, University of Pennsylvania, 209 South 33rd Street 
Philadelphia, PA 19104}
\email{blockj@math.upenn.edu}

\author{Nigel Higson}
\address{Department of Mathematics, Penn State University, University Park, PA 16802}
\email{higson@psu.edu}
\thanks{The second and third authors were partially supported by NSF grant DMS-1952669.}

\author{Jesus Sanchez Jr}
\address{Department of Mathematics, Penn State University, University Park, PA 16802}
\curraddr{}
\email{jxs1504@psu.edu}

\date{\today}



\begin{abstract}
We shall present a simplified version of a construction due to Denis Perrot that recovers the Todd class of the complexified tangent bundle of a smooth manifold from a JLO-type cyclic cocycle. The construction takes place within an algebraic framework, rather than the customary functional-analytic framework for the JLO theory.  The   series expansion for the exponential function is used in place of the heat kernel from the functional-analytic theory;  the Dirac operator chosen is far from elliptic; and a remarkable new trace discovered by Perrot replaces the operator trace.  In its full form Perrot's theory constitutes a wholly new approach to index theory. The account presented here here covers most but not all of his approach.  
\end{abstract}

\maketitle

\section{Introduction}
\label{sec-introduction}
The purpose of this paper is to give an exposition of some   remarkable ideas about index theory in the framework of  cyclic cohomology that Denis Perrot introduced about a decade ago \cite{Perrot2013}. Despite the originality of Perrot's work, his ideas  seem not yet to have  been studied in much detail, or at any rate not by many.  Our aim is to address this circumstance by  presenting a streamlined account  of many of the main ideas to as wide an audience as possible.

In \cite{Perrot2013}, Perrot gives a more or less  complete  account of the Atiyah-Singer index theorem using his new approach.  We shall not go that far.  Our aim instead will be to present Perrot's striking construction of the Todd class within the context of cyclic cohomology theory. In the course of doing so we shall come into contact with many of Perrot's most fascinating discoveries.

In later works Perrot has applied his new approach to index theory to new contexts---involving groupoids \cite{Perrot13b,Perrot16} and, in joint work with Rudy Rodsphon, foliations \cite{PerrotRodsphon14}.  The latter work settles in almost complete generality a longstanding open problem of Alain Connes and Henri Moscovici \cite{ConnesMoscovici95,ConnesMoscovici98}.  These applications  give ample evidence of the power of the new approach.  But once again  our aims will be  more limited, and we shall not discuss these developments here.

 Let $M$ be a smooth, closed manifold and let $\nabla$ be a torsion-free connection on $M$.  The curvature of $\nabla$   is a $2$-form $R$ on $M$ with values in $\End (TM)$, and the Todd form associated to $R$ is 
 \[
 \Todd (R) = \det \left (  \frac{R}{\exp(R)- 1}\right ) .
 \]
 This is a closed differential form on $M$ of mixed even degree.  
 Although its origins are in algebraic geometry, the Todd class is probably most famous for the role it plays in the Atiyah-Singer index theorem.  For instance, if $P$ is an elliptic $N{\times}N$ system of pseudodifferential operators on $M$, then
 \[
 \Index (P ) = \int _{S^*M} \operatorname{ch} (\sigma ) \Todd \bigl (  {R }\big / {2 \pi i}\bigr ) ,
 \]
where  $\sigma$ is the symbol of $P$, which is a smooth map from the cosphere bundle $S^*M$ into invertible $N{\times}N$ matrices, and $\operatorname{ch}(\sigma)$ is its Chern character, which is a closed differential form on $S^*M$ of mixed odd degree. The (normalized) Todd form is pulled back from $M$ to   $S^*M$. 
 
The purpose of this paper is to present, following Perrot, a new  construction of the Todd class using Dirac operators and cyclic cohomology.  
 
Let us recall very briefly  that if $A$ is a complex associative algebra, and if 
 \[
 C^p (A) = \bigl \{ \,\text{complex $(p{+}1)$-multilinear functionals on $A$}\,\bigr \} ,
 \]
 then there are differentials 
 \[
 b\colon C^p (A) \to C^{p+1}(A) 
 \quad \text{and} \quad 
  B\colon C^p (A) \to C^{p-1}(A) 
  \]
  with $(b{+}B)^2 = 0$.  A \emph{periodic cyclic cocycle} is a finitely supported family of $\{ \varphi _p\}_{p\ge 0}$  with $\varphi_p\in C^p(A)$ and 
  \[
  B\varphi_{p+1} + b\varphi_{p-1}=0\qquad  \forall p\ge 0.
  \]
If $A$ is the algebra of smooth functions on a smooth, closed  manifold, and if  $C$ is a  de Rham current in degree $p$, then the formula 
  \[
  \varphi _C (a^0,\dots, a^p) = \frac {1}{p!}\int _C a^0 da^1 \cdots da^p
  \]
  defines an element  $\varphi_C\in C^p(A) $ with 
  \[
  b \varphi _C = 0 \quad \text{and} \quad B \varphi _C = \varphi_{d'C} ,
  \]
  where $d'$ is the de Rham differential on currents. So if $C$ is closed, then $\varphi_C$ is a periodic cyclic cocycle (concentrated in degree $p$).   
  
It follows from the above that  if we  pull the Todd form back from $M$ to $S^*M$, as in the index theorem, and if we denote by $\Todd^{2q} (R)$  the degree $2q$-component, then for $p{+}2q = 2\dim (M) {-}1$
the formula 
\[
\varphi_p (a^0,\dots, a^p ) = \frac 1 {p!}\int _{S^*M} a^0 da^1 \cdots da^{p} \Todd^{2q} (R) \qquad (\text{$p $ odd}),
\]
defines a (mixed degree) periodic cyclic cocycle.
  
It is this cocycle that Perrot constructs in a new way. From the point of view of the formalism of cyclic theory, Perrot's construction is a variation on the famous JLO formula
 \[
 \STr \Bigl ( \int _{\Sigma^p} a^0 \exp (-s_0 \Dirac^2 ) [\Dirac, a^1] \exp (-s_1 \Dirac^2 )
 \cdots  [\Dirac,a^p]  \exp (-s_p \Dirac^2 )\, ds 
 \Bigr ) .
 \]
See \cite{JLO88,Quillen88}. What makes it noteworthy  are:
\begin{enumerate}[\rm (i)]

\item   Perrot's definition of the Dirac operator --- among other things it is \emph{not} an elliptic  operator; 

\item   Perrot's approach to the heat operators $\exp (-s \Dirac ^2)$ --- they are defined \emph{algebraically}, using the series expansion for the exponential function; and 

\item Perrot's definition of the (super)trace in the JLO formula --- it is not by any means an operator trace, or even a residue trace  \cite{Wodzkicki87,Kassel89}, although it is related to the latter.
\end{enumerate}
Our goal in this paper is to   discuss each of these points in detail.  But to conclude our introduction, let us indicate  the parts of Perrot's work that we shall \emph{not} cover in the paper.
 
In \cite{Perrot2013}, Perrot  works throughout with the algebra of order zero classical pseudodifferential symbols on a smooth, closed manifold $M$.  Our more limited goals in this paper make it possible for us  to substantially simplify matters by working instead with the \emph{commutative} associated graded algebra of polyhomogeneous smooth functions on the cotangent bundle.  We shall be able to explain Perrot's  realization of the Todd form while working exclusively in this context, thereby replicating a major part of  \cite{Perrot13b}. But to complete the proof of the index theorem it would be necessary to lift the entire discussion from polyhomogeneous functions to symbols. 
 
In the context of pseudodifferential symbols, Perrot actually constructs \emph{two} periodic cyclic cocycles. One is effectively the cocycle that we shall construct and compute in this paper.  The other is constructed in a very similar fashion, but Perrot shows that it is precisely the so-called Radul cocyle \cite{Radul91}, which produces the analytic index \cite{Nistor97,Perrot12}. He also proves  that   his two cocycles are cohomologous. So, putting everything together, the  computations of the two cocycles amount to a proof of the Atiyah-Singer index theorem for elliptic pseudodifferential systems.
 
\paragraph{\emph{A note on the text.}}  Nearly all the arguments presented below are adapted from Perrot's paper \cite{Perrot2013}, some of them very closely.  In places, for example in Section~\ref{sec-dirac-operator} where we present a  simple global construction of the Dirac operator, we have been able to take advantage of our simplified context to streamline some of Perrot's constructions.  In other places, for example in  our treatment of the coordinate independence of Perrot's trace in Section~\ref{sec-perrot-trace}, we have chosen an   approach that deviates from \cite{Perrot2013}, typically because we felt it was illuminating to do so.  But overall we owe a huge debt to \cite{Perrot2013}. 

\paragraph{\emph{Acknowledgement.}}  Beyond acknowledging our indebtedness to Perrot, it is a pleasure to thank Rudy Rodsphon for sharing his insights  into  Perrot's work    during a number of stimulating conversations.

\section{Perrot's Dirac Operator} 
\label{sec-dirac-operator}
Throughout the paper,     $M$ will always be  a smooth manifold without boundary.  We shall always denote by $n$ the dimension of $M$.  In various places (where we need to integrate over $M$) we shall in addition assume that $M$ is compact.  We shall be working with vector fields, differential forms, etc, on $M$; all these will  be assumed to be \emph{smooth}.

\begin{definition} 
Throughout the paper, we shall denote by $\Tcirc M$ the complement of the set of zero covectors in the cotangent bundle $T^*M$: 
\[
\Tcirc M = T^*M \setminus M.
\]
\end{definition}

In this section we shall work with a fixed   affine connection $\nabla$ on the tangent bundle of $M$ (eventually we will require $\nabla$ to be torsion-free). Following Perrot (but with some variations, as explained in the introduction) we  are going to use $\nabla$ to construct  a Dirac operator on the total space $\Tcirc M$ that acts on sections of the pullback to $\Tcirc M$ of the exterior algebra bundle of $M$. But we should say at the outset that this Dirac operator will \emph{not} be elliptic, and in particular it will not be a typical Dirac operator from index theory.

The construction  requires  some facts about  horizontal vector fields on $\Tcirc M$.  In what follows we shall denote by    $\pi\colon \Tcirc M\to M$  the standard projection mapping.

\begin{definition}
\label{def-of-alpha}
We shall identify
vector fields on $M$ with fiberwise linear smooth functions on $T^*M$ via the isomorphism of $C^\infty (M)$-modules 
\[
\alpha \colon \bigl\{ \, \text{vector fields on $M$} \,\bigr \} 
\stackrel \cong \longrightarrow 
\Bigl\{ \, \parbox{1.8in}{\begin{center}fiberwise linear smooth functions  on $T^*M$\end{center}} \,\Bigr \} 
\]
defined by 
$\alpha (X)(\omega) = \omega (X)$ for every $1$-form $\omega$ on $M$.
In local coordinates,
\[
\alpha \Bigl ( \frac{\partial}{\partial x^i}\Bigr ) = \xi_i ,
\]
where $\xi_1,\dots , \xi_n$ are the usual fiberwise (linear) coordinate functions on $T^*M$ dual to $x^1,\dots, x^n$, defined by $\xi_i (\alpha) = \alpha (\partial/\partial x^i)$.
\end{definition}

\begin{definition}
Let $X$ be a  vector field on $M$. Its \emph{horizontal lift} to the total space of the cotangent bundle $T^*M$ (with respect to the connection $\nabla$) is the vector field $X^H$ on $T^*M$ that is $\pi$-related to $X$ and that is characterized by the further requirement that 
\[
X^H (\alpha(Y)) = \alpha(\nabla_X Y)  ,
\]
for every vector field $Y$ on $M$.
\end{definition} 

In local coordinates, if we introduce the usual Christoffel symbols by
\[
\nabla_{{\partial}/{\partial x^i}} \frac{\partial}{\partial x^j}
= \Gamma_{ij}^k\, \frac{\partial}{\partial x^k} 
\]
(with the usual summation convention), then 
\begin{equation}
    \label{eq-local-formula-for-horizontal-vector-field}
  \Bigl (\frac{\partial}{\partial x^i}\Bigr )^H 
=  \frac{\partial}{\partial x^i} +  \Gamma_{ij}^k \,\xi_k \frac{\partial}{\partial \xi_j}, 
\end{equation}
where $\xi_1,\dots , \xi_n$ are, as above, the fiberwise linear coordinate functions on $T^*M$ dual to $x^1,\dots, x^n$.

\begin{definition}
\label{def-of-gamma}
We introduce the following additional identification:
we define an isomorphism of $C^\infty (M)$-modules 
\[
\gamma \colon \operatorname{End}\bigl ( TM\bigr ) 
\stackrel \cong \longrightarrow 
\Bigl \{ \,
\text{vertical vector fields $a_{ij} \,\xi _i  \,{\partial} /{\partial \xi_j}$ on $T^*M$}
\,\Bigr \} 
\]
by the formula 
$\gamma (A)(\alpha (X)) = \alpha (A(X))$. 
Thus in local coordinates,  if $A \bigl ( {\partial} /{\partial x^i}\bigr ) = a_{ij} {\partial} /{\partial x^j}$,
then 
\[
\gamma (A) = a_{ij}\, \xi_i\, \frac{\partial} {\partial \xi _j} .
\]
\end{definition}

\begin{lemma}
 \label{lem-horiz-vector-fields-and-curvature}
 If $X$ and $Y$ are vector fields on $M$, then  
\[
[X^H,Y^H]-[X,Y]^H = \gamma \bigl (R(X,Y)\bigr )
\]
\end{lemma}

\begin{proof}
Both sides of the identity are vertical vector fields, so it suffices to check that both sides agree on the functions $\alpha(Z)$ in Definition~\ref{def-of-alpha}.  The identity in this case is an immediate consequence of the definitions.
\end{proof}

\begin{definition} 
Form the exterior algebra bundle $\Lambda^* T^*M$ over $M$. Equip it with the   connection induced from $\nabla$, and then pull  the bundle and its connection back from $M$ to the manifold $\Tcirc M$. For brevity, we shall write 
\[
S = \pi ^* \Lambda ^* T^*M
\]
for the bundle, and keep the  notation $\nabla$ for the connection.
\end{definition}

In order to express  the connection on $S$ in local coordinates it will be convenient to introduce the following notation: 
\begin{definition}
Given local coordinates $x^1,\dots, x^n$   on an open set $U\subseteq M$, define operators $\psi^i$ and $\overline \psi_j$ acting on sections of $S$  by
\[
\begin{aligned}
\psi^i & = \text{left exterior multiplication by $dx^i$ on $\pi^* \Lambda^* T^*M$} \\
\overline \psi _j & = \text{contraction by $\partial / \partial x^j$ on $\pi^* \Lambda^* T^*M$}.
\end{aligned} 
\]
\end{definition}

Using  the above notation, along with the previously introduced Christoffel symbols $\Gamma _{ij}^k$, the induced connection on $\pi^* \Lambda ^* T^*M $ is 
\begin{equation}
    \label{eq-formula-for-induced-connection}
\begin{aligned}
\nabla _{(\partial/\partial x^i)^H} & = \frac{\partial}{\partial x^i} + \Gamma_{i\ell}^k\, \xi_k\, \frac{\partial }{ \partial \xi _\ell} - \Gamma_{i\ell}^k \,\psi^\ell\, \overline \psi_k
\\
\nabla _{\partial/\partial \xi_j } & = \frac{\partial}{\partial \xi_j} .
\end{aligned}
\end{equation}

\begin{definition}
\label{def-d-vert-and-d-horiz}
Given local coordinates on an open set $U\subseteq M$, define operators 
\[
\Dirac_{\mathrm{horiz}}, \Dirac_{\mathrm{vert}}\colon \Gamma \bigl (T^*U, \pi^* \Lambda ^* T^*M\bigr )
\longrightarrow 
\Gamma \bigl (T^*U, \pi^* \Lambda ^* T^*M\bigr )
\]
using the   formulas
\[
\Dirac_{\mathrm{horiz}} =   \psi^i \nabla _{(\partial/\partial x^i)^H}\quad \text{and} \quad 
\Dirac_{\mathrm{vert}} =   \overline \psi_j \, \nabla_{\partial/\partial \xi_j} .
\]
\end{definition}

Using the Christoffel symbols for $\nabla$ and the standard local frame for $S$ associated to a coordinate system, we find that 
\[
\Dirac_{\mathrm{horiz}} 
     =  \psi^i  \frac{\partial}{\partial{x^i}}+ \Gamma^k_{ij} \psi ^i\xi_k\frac{\partial}{\partial{\xi_j}}-\Gamma^k_{ij}\psi^i\psi^j\overline \psi_k  
     \quad \text{and} \quad 
     \Dirac_{\mathrm{vert}} = \overline \psi^j \frac{\partial}{\partial \xi_j}.
\]

\begin{lemma} 
The operators 
$\Dirac_{\mathrm{horiz}}$ and  $\Dirac_{\mathrm{vert}}$ 
are independent of the choice of local coordinate system used to define them.
\end{lemma}

\begin{proof}
This follows from the $C^\infty(M)$-linearity of the map  $\alpha$.
\end{proof}

\begin{lemma}
\label{lem-vanishing-term-for-torsion-free-connection}
If the connection $\nabla$ on $TM$ is torsion-free, then  
\[
\Dirac_{\mathrm{horiz}} 
      =   \psi^i  \frac{\partial}{\partial{x^i}}+ \Gamma^k_{ij} \psi ^i\xi_k\frac{\partial}{\partial{\xi_j}} 
\]
in any local coordinate system.
\end{lemma}

\begin{proof}
This is because $\Gamma_{ij}^k = \Gamma^k_{ji}$ for a torsion-free connection, so that the third term in the defining formula 
\[
\Dirac_{\mathrm{horiz}} 
     =  \psi^i  \frac{\partial}{\partial{x^i}}+\Gamma^k_{ij}\psi ^i\xi_k\frac{\partial}{\partial{\xi_j}}-\Gamma^k_{ij}\psi^i\psi^j\overline \psi_k  
\]
is zero in this case.
\end{proof}

\begin{lemma}
\label{lem-square-of-d-horiz}
If the connection $\nabla$ on $TM$ is torsion-free, then
\[
\Dirac_{\mathrm{horiz}}^2    =   \tfrac{1}{2}\psi^i\psi^j \gamma \bigl (  R  ({\partial}/{\partial x^i}, {\partial}/{\partial x^j}  )\bigr )
\]
 in any local coordinate system.
 \end{lemma}
 
 \begin{proof}
 Fix a local coordinate system on $U\subseteq M$.  Give $S$  the local frame  associated to these coordinates, and define operators 
 \[
 \partial^H_i = \frac{\partial}{\partial x^i} + \Gamma^k_{i\ell} \xi_k \frac{\partial}{\partial \xi_\ell}
 \]
 on smooth sections of $S$ (these  are not the operators  $\nabla_{(\partial/\partial x^i)^H}$, nor  do they have any particularly special meaning; they are merely introduced for the sake of the computation).  We now calculate from Lemma~\ref{lem-vanishing-term-for-torsion-free-connection} that 
 \[
 \begin{aligned}
 \Dirac_{\mathrm{horiz}}^2 
    & =   \psi ^i \partial^H_i   \cdot   \psi ^j \partial^H_j   \\
    & = \psi^i \psi^j \partial^H_i\partial^H_j \\
    & = \sum_{i<j} \psi^i \psi^j \bigl ( \partial^H_i\partial^H_j - \partial^H_j\partial^H_i \bigr ) . \\
 \end{aligned}
 \]
 But it follows from Lemma~\ref{lem-horiz-vector-fields-and-curvature} that 
 \[
 \partial^H_i\partial^H_j - \partial^H_j\partial^H_i = \gamma \bigl (R (\partial/\partial x^i, \partial/\partial x^j)\bigr ),
 \]
 and as a result 
 \[
  \begin{aligned}
  \Dirac_{\mathrm{horiz}}^2 
    & =  \sum_{i<j} \psi^i \psi^j \gamma \bigl (R (\partial/\partial x^i, \partial/\partial x^j)\bigr ) \\
    & = \tfrac 12 \psi^i \psi^j \gamma \bigl (R (\partial/\partial x^i, \partial/\partial x^j)\bigr ) , \\
  \end{aligned}
 \]
 as required.
 \end{proof}

 \begin{definition}
 \label{def-gamma-of-r}
We shall write 
\[
\gamma(R) =  \tfrac{1}{2}\psi^i\psi^j \gamma \bigl (  R  ({\partial}/{\partial x^i}, {\partial}/{\partial x^j} )\bigr ) .
\]
This is   a first-order linear differential operator acting on section of $S$,  independent of the choice of coordinates. 
More explicitly, if we write 
\[
R \Bigl(\frac{\partial}{\partial x^i}, \frac{\partial}{\partial x^j}\Bigr )\, \frac{\partial}{\partial x^\ell} =   R^k_{ij\ell}\, \frac {\partial}{\partial x^k} ,
\]
then 
\[
\gamma(R)   =   \frac{1}{2}\psi^i\psi^jR_{ij\ell}^k \xi_k \frac{\partial}{\partial \xi _\ell} .
\] 
\end{definition}

In this paper, following Perrot  we shall be  concerned with  the ``Dirac operator'' $D=\Dirac_{\mathrm{horiz}}{+}\Dirac_{\mathrm{vert}}$ (or actually a small variation of this; see Definition~\ref{def-dirac-operator}) and its square. Lemma~\ref{lem-square-of-d-horiz} is obviously relevant to the computation of the square, as is the simple formula 
\begin{equation}
\label{eq-square-of-vertical-operator}
    D^2_{\mathrm{vert}} = 0. 
\end{equation}
What remains is to compute the cross-term 
\[
\{\Dirac_{\mathrm{horiz}} , \Dirac_{\mathrm{vert}}\} =  \Dirac_{\mathrm{horiz}} \Dirac_{\mathrm{vert}} {+} \Dirac_{\mathrm{vert}}\Dirac_{\mathrm{horiz}}.
\]

 \begin{lemma}
 \label{lem-local-formula-for-cross-term}
If the connection $\nabla$ on $TM$ is torsion-free, then 
\[
\{\Dirac_{\mathrm{horiz}} , \Dirac_{\mathrm{vert}}\}
= 
\nabla_{\partial_{\xi_i}}\nabla_{ ( \partial _{x^i})^H}  
\]
\textup{(}with the summation convention, as always\textup{)} in any local coordinate system.
\end{lemma}

\begin{proof}
We shall use the same operators $\partial_i^H$ introduced in the proof of Lemma~\ref{lem-square-of-d-horiz}.  With these, we may write 
\begin{equation}
\label{eq-cross-term-1}
\begin{aligned}
\{\Dirac_{\mathrm{horiz}} , \Dirac_{\mathrm{vert}}\}
    & = \bigl ( \psi ^i \partial_i ^H\bigr ) \bigl (\overline \psi_j \frac{\partial}{\partial \xi _j} \bigr )
    +
    \bigl (\overline \psi_j \frac{\partial}{\partial \xi _j} \bigr )\bigl ( \psi ^i \partial_i ^H\bigr ) \\
    & = \psi ^i  \overline \psi_j \partial_i ^H  \frac{\partial}{\partial \xi _j}  
    +
     \overline \psi_j  \psi ^i  \frac{\partial}{\partial \xi _j}   \partial_i ^H   .
\end{aligned} 
\end{equation}
Next, we compute that 
\[
    \partial_i ^H  \frac{\partial}{\partial \xi _j}
      =  \Bigl ( \frac{\partial}{\partial x^i} + \Gamma^k_{i\ell} \xi_k \frac{\partial}{\partial \xi_\ell}\Bigr )  \frac{\partial}{\partial \xi _j}
     =
     \frac{\partial}{\partial \xi _j}   \partial_i ^H 
    - \Gamma^j_{i\ell}   \frac{\partial}{\partial \xi_\ell} ,
\]
and as a result we find from \eqref{eq-cross-term-1} that 
\begin{equation}
    \label{eq-cross-term-2}
\begin{aligned}
\{\Dirac_{\mathrm{horiz}} , \Dirac_{\mathrm{vert}}\}
    & = \bigl (  \psi ^i  \overline \psi_j{+}  \overline \psi_j \psi ^i \bigr)   \frac{\partial}{\partial \xi _j}  \partial_i ^H  -  \Gamma^j_{i\ell} \psi ^i    \overline \psi_j  \frac{\partial}{\partial \xi_\ell} 
    \\
    & =  \frac{\partial}{\partial \xi _i} 
    \Bigl ( \frac{\partial}{\partial x^i} + \Gamma^k_{i\ell} \xi_k \frac{\partial}{\partial \xi_\ell}\Bigr )
    -  \Gamma^j_{i\ell} \psi ^i   \overline \psi_j  \frac{\partial}{\partial \xi_\ell} 
    \\
    & =  \frac{\partial^2\,\,\,}{\partial \xi_i\partial x^i} + \Gamma^k_{i\ell} \xi_k \frac{\partial^2\,\,\,}{\partial \xi_i \partial \xi_\ell} 
    + \Gamma^i_{i\ell}   \frac{\partial}{\partial \xi_\ell}
    - \Gamma^j_{i\ell}  \psi ^i  \overline \psi_j  \frac{\partial}{\partial \xi_\ell} 
    .
\end{aligned}
\end{equation}
On the other hand, 
\begin{equation}
\label{eq-cross-term-3}
\begin{aligned}
\nabla_{\partial _{\xi _i}} \nabla _{(\partial/\partial x^i)^H} 
    & = \frac{\partial}{\partial \xi_i} \left ( 
    \frac{\partial}{\partial x^i} + \Gamma_{i\ell}^k\, \xi_k\, \frac{ \partial}{\partial \xi _\ell} - \Gamma_{i\ell}^k \,\psi^\ell\, \overline \psi_k
    \right )\\
    & =   
    \frac{\partial^2 }{\partial \xi_i\partial x^i} + \Gamma_{i\ell}^k\, \xi_k\, \frac{ \partial^2}{\partial \xi _i\partial \xi _\ell}
     + \Gamma_{i\ell}^i\,   \frac{ \partial}{ \partial \xi _\ell}
     - \Gamma_{i\ell}^k \,\psi^\ell\, \overline \psi_k
    \frac{\partial}{\partial \xi_i} . \\
\end{aligned}
\end{equation}
Since $\nabla$ is torsion-free, $\Gamma^k_{i\ell} = \Gamma^k_{\ell i}$, and so we find that \eqref{eq-cross-term-2} and \eqref{eq-cross-term-3} are equal, as required.
 \end{proof}

If we  write 
\[
\nabla^2 = \nabla_{\partial_{\xi_i}}\nabla_{ (\partial _{x^i})^H}
\]
(this is coordinate-independent), then the computations in this section may be summarized as follows:  
\begin{equation}
\label{eq-summary-of-Dirac-squared-formulas}
  \Dirac_{\mathrm{horiz}} ^2 = \gamma (R),\quad  \Dirac_{\mathrm{vert}} ^2 = 0
    \quad\text{and}\quad 
    \{\Dirac_{\mathrm{horiz}} , \Dirac_{\mathrm{vert}}\}
    = \nabla^2 .
\end{equation}
These identities will be crucial in the sequel. 

\begin{remark}
Using the connection $\nabla$  we may identify   the tangent bundle of $\Tcirc M$ with the pullback of $TM {\oplus} T^*M$ to $\Tcirc M$. If we equip the latter with its canonical symmetric bilinear form, which of course is nondegenerate but of indefinite signature,  then $S$ carries an irreducible representation of the associated bundle of Clifford algebras, and the operator $\Dirac_{\mathrm{horiz}} {+} \Dirac_{\mathrm{vert}}$ is a Dirac operator in the usual sense, although since the bilinear form is not definite, the operator is not elliptic.
\end{remark}

\section{Infinite-Order Polyhomogeneous Operators}

In this section we shall describe the infinite-order linear partial differential operators (presented as formal  series) with which we shall work in the rest of the paper.  We shall begin by considering  differential operators that act on scalar functions, but then we shall turn  to differential operators acting on spinors, which will be our principal objects of interest.

\begin{definition}
\label{def-polyhomogeneous-function}
We shall denote by $E$ the Euler vector field on $\Tcirc M$, so that 
\[
(Ef )(\alpha) = \frac{d}{dt}\Big \vert _{t  =0} f (e^t \alpha) 
\]
for all smooth functions $f$ on $\Tcirc M$ and all $\alpha \in \Tcirc M$.  We shall say that a smooth function on $\Tcirc M$ is \emph{homogeneous of degree $k\in \Z$}  if it is a $k$-eigenfunction for the action of $E$, or in other words if 
\[
f (e^t \alpha) = e^{kt} f (\alpha)
\]
for all $t\in \R$ and all $\alpha \in \Tcirc M$. We shall say that $f$ is 
 \emph{polyhomogenous} if it is a (finite) linear combination of smooth homogeneous functions of various integer degrees.  We shall write 
 \[
\Poly^k (M) = \bigl \{ \, \text{Smooth degree-$k$ homogeneous functions on $\Tcirc M$}\,\bigr  \} 
\]
and
\[
\Poly (M) = \bigl \{ \, \text{Smooth polyhomogeneous functions on $\Tcirc M$}\,\bigr  \} .
\]
\end{definition}

\begin{definition}
\label{def-pdo-algebra}
A linear partial differential operator $D$ on $\Tcirc M$,  acting on smooth functions on $\Tcirc M$, is \emph{homogeneous of degree $k\in \Z$} if 
 $\ad_{E} (D) = k D$, where $E$ is the Euler vector field, and \emph{polyhomogeneous} if it is a (finite) linear combination of homogeneous operators of various integer degrees.  
We shall denote by $\PDO(M)$ the algebra of polyhomogeneous linear partial differential operators on  $\Tcirc M$.
\end{definition}

\begin{remark}
\label{rem-standard-form-of-pdo}
In local coordinates,  the  polyhomogeneous operators   are linear combinations of operators of the form
\[
  f_{\alpha \beta } \, \frac{\partial^\alpha}{\partial x^\alpha} \, \frac{\partial ^\beta}{\partial \xi_\beta} ,
\]
with $f$ a  polyhomogeneous smooth function on $\Tcirc M$.
\end{remark}

\begin{definition}
We shall denote  by  $F$ the  field of  formal complex Laurent series $\sum _{-N}^\infty a_k \varepsilon ^k$ (with only  finite singular parts). We shall denote by 
$\FPDO(M)$ 
the algebra of formal Laurent series in $\varepsilon$  with coefficients from $\PDO(M)$ (again, with only  finite singular parts). 
\end{definition}

We turn now to one of  the main constructions of the paper, which is that of  a subspace 
\[
\FB(M)\subseteq \FPDO(M)
\]
on which will be defined a crucial trace functional. The subspace  is not a subalgebra, let alone an ideal, but it is for instance a bimodule under the left and right actions of $\PDO(M)$ on $\FPDO(M)$ (the name is supposed to suggest this, and Perrot calls his version, on which ours is very closely modelled, the \emph{bimodule of trace class operators}).  We shall define the bimodule in this section and construct the trace in the next section.

The space $\FB(M)$ will be defined in stages, and to begin we shall work in a fixed coordinate system on an open subset $U\subseteq M$. We shall use the following operator extensively:

\begin{definition}
\label{def-laplacian-in-local-coordinates}
Let $U$ be an open subset of $M$ and let $x^1,\dots, x^n$ be coordinates on $U$.  Given these coordinates, the associated \emph{Laplacian} on $\Tcirc U$ is
\[
\Delta = \frac{\partial^2}{\partial x^i\partial\xi_i}
\]
(summation convention).
\end{definition}

\begin{remarks} The operator $\Delta$ is \emph{not} invariant under changes of coordinates (except for linear changes of coordinates).  In addition, the operator $\Delta$ is \emph{not} elliptic.
\end{remarks}

Following Perrot  we   introduce the following increasing filtration on the algebra $\PDO (U)$.

\begin{proposition/definition}
\label{prop-def-bimodule-order-scalar-case}
Let $U$ be an open subset of  $M$ and  let $x^1,\dots, x^n$ be coordinates on $U$.  Define an increasing filtration of the algebra $\PDO(U)$, indexed by $\Z$, and an associated notion of order  on partial differential operators on $U$, that we shall call \emph{bimodule order}, as follows. The bimodule filtration is constructed by  considering at first all increasing algebra filtrations 
\[
\cdots \subseteq \PDO_p (U) \subseteq \PDO_{p+1}(U) \subseteq \cdots \subseteq \PDO(U)
\]
indexed by $p\in \Z$ and with $\partial/\partial x^i\in \PDO_3(U)$, $\partial/\partial \xi_i\in \PDO_{-1}(U)$, and $f\in \PDO_{2p} $ if $f$ is homogeneous of degree $p$.  Then define the $p$th bimodule filtration space to be the intersection of all $\PDO_p (U)$ over all these filtrations. Of course we  write $\bimodorder (T)\le p$ if $T$ belongs to this intersection.  We have the following identities:
\[
\begin{aligned}
 \bimodorder( \partial / \partial x^i )  & = 3 \qquad \,\,\,\forall i\\
 \bimodorder( \partial / \partial \xi_j )  &  = -1\qquad \forall j \\ 
\bimodorder (f)  & =  2p \qquad \,\, \text{if $\operatorname{degree}(f) = p$.} 
\end{aligned}
\]
\end{proposition/definition}

\begin{proof} A straightforward computation in local coordinates.
\end{proof}

\begin{remark}
We shall introduce a second,  quite different filtration on differential operators in Section~\ref{sec-perrot-order}, and it should not be confused with the  bimodule filtration defined here.  
\end{remark}

\begin{definition} 
\label{def-definition-of-fbu}
We shall denote by $\FB (U) \subseteq \FPDO(U)$ the space of all Laurent series 
\[
\sum _{k} \varepsilon ^ k D_k  \exp (\varepsilon \Delta)  \qquad (D_k \in \PDO (U)) 
\]
for which there exists $N>0$ with 
\[
 \bimodorder(D_k) \le k + N  ,
\]
for all $k$.
\end{definition} 

\begin{remark}
It is  not so easy to explain at the outset the reason for defining the bimodule filtration, and hence $\FB(U)$, in the way that we just have.  In fact there is a certain amount of flexibility in how one might define $\FB(U)$, which needs to be small enough that the trace functional to be discussed in the next section is well-defined, and large enough that it is closed under some simple operations, most notably conjugation by $\exp(\varepsilon \Delta)$.  The definition given seems to be the simplest one that meets these requirements.\footnote{One might perhaps compare this to the situation with traceable Hilbert space operators in index theory, where one might work with smoothing operators, or trace-class operators, or  operators whose singular value sequence has rapid decay, and so on.}  
\end{remark}

The following is another easy computation:

\begin{lemma} 
\label{lem-bimodule-order-and-ad-delta}
If $T$ is any operator in $\PDO(U)$, then 
\[
\pushQED{\qed} 
\bimodorder \bigl ([\Delta, T]\bigr ) \le \bimodorder (T) + 1.
\qedhere
\popQED
\]
\end{lemma}

Now suppose that $T\in FPDO(U)$ and that  $s\in \R$. The conjugated operator  
\[
\operatorname{Ad}_{\exp (\varepsilon s \Delta}(T) = \exp (\varepsilon s \Delta) T \exp (-\varepsilon s \Delta)
\in \FPDO (U) 
\]
may be alternatively expressed as 
\[
\operatorname{Ad}_{\exp (\varepsilon s \Delta )}(T) 
=
\exp (\varepsilon s \operatorname{ad}_{ \Delta})(T) 
=
\sum _{k=0}^\infty \frac{s^k \varepsilon ^k}{k!} \ad_{\Delta}^k (T)
\]
(thanks to the powers $\varepsilon ^k$,  the sum makes sense in $\FPDO(U)$; to use language that will be introduced in the next section, the sum  is convergent in the adic topology). Using Lemma~\ref{lem-bimodule-order-and-ad-delta} we find that: 

\begin{corollary} 
\label{cor-closure-under-ad-exp-delta}
If $T\in \PDO(U)$, then for every $s\in \R$  we may write 
\[ 
\operatorname{Ad}_{\exp (\varepsilon s \Delta )} ( T) =\sum _{k=0}^\infty \varepsilon^k T_k
\]
where $T_k\in \PDO(U)$ and  $\order (T_k) \le k + \order(T)$.  In particular the subspace $\FB(M)\subseteq \FPDO (M)$ is closed under left  and right multiplication by elements of $\PDO(M)$ in the algebra  $\FPDO(U)$. \qed
\end{corollary}

\begin{corollary} 
\label{cor-closure-under-ad-exp-s-delta}
If $X\in \FB(U)$, then   $\operatorname{Ad}_{\exp (\varepsilon s \Delta )} ( X)  \in \FB(U)$  for every $s\in \R$.    \qed
\end{corollary}

The bimodule filtration on $\PDO(U)$  is independent of the choice of coordinates on $U$; this can be checked by the usual change-of-coordinates formulas.   The following result is a bit more involved:

\begin{lemma}
\label{lem-exp-delta-prime-in-fb}
Let $U$ be an open subset of $M$,  let $x^1,\dots, x^n$ be coordinates on $U$, and let $\Delta $ and $\FB(U)$ be the associated Laplacian and bimodule.  If $\Delta'$ is the Laplacian associated to another coordinate system on $U$, then \[
\exp(\varepsilon \Delta') \in \FB(U).
\]
\end{lemma}

\begin{proof} 
The Laplacian $\Delta'$ in a new coordinate system has the form
\[
\Delta' = \Delta + a_{i} \frac{\partial}{\partial \xi_i} + b_{jk\ell} \xi_j \frac{\partial^2}{\partial \xi_k\partial \xi_\ell}
\]
in terms of the original coordinate system, where $a_i$ and $b_{jk\ell}$ are smooth functions on $U$, and so 
$\bimodorder (\Delta {-}  \Delta') \le 0$.

To proceed, within $\FPDO (U)$ there is a Duhamel-type formula 
\[
\exp (\varepsilon \Delta ') - \exp (\varepsilon \Delta) 
=
\varepsilon \int _0^1 \exp (\varepsilon s \Delta ') (\Delta'{-}\Delta)\exp (\varepsilon (1{-}s) \Delta)\, ds.
\]
The integrand, viewed as a power series in $\varepsilon$, is a polynomial function of $s$ in each degree, with values in $\PDO(U)$. So the definition of the integral presents no problems (and the formula may be proved in the usual way).  

By iterating the formula we obtain the perturbation series 
\begin{multline*}
\exp (\varepsilon \Delta ' ) =
\sum_{p=0}^\infty \varepsilon ^p \int _{\Sigma^p}
\exp (\varepsilon s_0 \Delta) R\exp (\varepsilon s_1\Delta)R\cdots 
\\
\cdots \exp(\varepsilon s_{p-1}\Delta)R \exp (\varepsilon s_p\Delta)\, ds ;
\end{multline*}
in $\FPDO(U)$, where $R= \Delta'{-}\Delta$ and where $\Sigma^p$ is the standard $p$-simplex.  Writing the above as 
\begin{multline*}
\exp (\varepsilon \Delta ' ) =
\sum_{p=0}^\infty \varepsilon ^p \int _{\Sigma^p}
\operatorname{Ad}_{\exp (\varepsilon s_0\Delta)} (R)
\operatorname{Ad}_{\exp (\varepsilon (s_0{+}s_1)\Delta)} (R)\cdots 
\\
\cdots \operatorname{Ad}_{\exp (\varepsilon (s_0{+}s_1{+}\cdots{+}s_{p-1})\Delta)} (R)
\, ds \cdot \exp (\varepsilon \Delta)
\end{multline*}
and using Corollary~\ref{cor-closure-under-ad-exp-delta}   we obtain the result.
\end{proof} 

\begin{theorem} 
The subspace $\FB(U)\subseteq \FPDO (U)$ is independent of the choice of local coordinates.
\end{theorem}

\begin{proof}
First, Lemma~\ref{lem-exp-delta-prime-in-fb}  shows that $\exp (\varepsilon \Delta') \in \FB(U) $.  Second, because the bimodule order is obtained from an algebra filtration, $\FB(U)$ is in fact closed under left multiplication by series $\sum_k \varepsilon ^k T_k$ with $T_k\in \PDO(U)$ and $\bimodorder(T_k) \le k + N  $  for some $N$ and all $k$.  These two observations show that 
the $\FB$-bimodule defined using $\Delta'$ is included in the $\FB$-bimodule defined using $\Delta$.  By symmetry, the two are in fact equal.
\end{proof}

We can now define a global version of the bimodule $\FB(U)$:

\begin{definition}
We denote by $\FB(M)\subseteq \FPDO(M)$ the subspace of all elements that restrict to elements in $\FB(U)$ whenever $U$ is an open subset of $M$ that supports a coordinate system.
\end{definition}

We shall now make the minor changes needed to  carry over  all of the results  above from operators that act on scalar functions to operators that act on sections of the bundle $S$.   The most direct way to do proceed is to adapt the notion of bimodule order to this context: 

\begin{definition}
\label{def-bimodule-order-spinor-case}
Let $U$ be an open subset of  $M$ and  let $x^1,\dots, x^n$ be coordinates on $U$. Trivialize the bundle $S$ over $U$ using the associated frame for the exterior algebra bundle over $U$, and  define an increasing filtration of the algebra $\PDO(U,S)$, indexed by $\Z$, and an associated notion of \emph{bimodule order} on partial differential operators on $U$, by
\[
\begin{aligned}
\bimodorder( \partial / \partial x^i )  & \le 3 \qquad \forall i\\
 \bimodorder( \partial / \partial \xi_j )  &  \le - 1\qquad \forall j \\ 
\bimodorder (f)  & \le  2p \qquad  \text{if $f$ is a scalar function of degree $p$} \\
\bimodorder(\psi^i) & \le 0 \qquad \forall  i\\
\bimodorder(\overline \psi_j) & \le 0  \qquad \forall j\\
\end{aligned}
\]
\end{definition}

Thus the bimodule order on $\PDO(U,S)$ does not involve the bundle $S$ at all.  It is coordinate-independent, and all of the results and definitions given above pass over to the context of operators on $S$ without any change at all (of course one must replace homogeneous scalar functions on $\Tcirc U$ with homogeneous endomorphisms of $S$, noting that since $S$ is pulled back from $M$, the notion of homogeneity makes sense here, as it does  for \emph{any} pull-back bundle, since the operation of scalar multiplication on the fibers of $\Tcirc M$ lifts canonically to an action on the bundle).

To conclude this section we shall situate the Driac operator and its exponential within the bimodule  $\FB(M,S)$.

\begin{definition}
\label{def-dirac-operator}
The \emph{Dirac operator} on $\Tcirc M$, associated to a torsion-free connection on $M$, is the operator 
\[
\Dirac = \varepsilon \Dirac _{\mathrm{horiz}} + \Dirac_{\mathrm{vert}} ,
\]
whose summands are described in Definition~\ref{def-d-vert-and-d-horiz}. 
\end{definition}

It follows from the formulas in \eqref{eq-summary-of-Dirac-squared-formulas} that 
\[
\Dirac ^2 = \varepsilon \nabla ^2 + \varepsilon ^2 \gamma (R),
\]
and therefore the exponentials 
\[
\exp (s \Dirac ^2 ) = \sum_{k=0}^\infty \frac{s^k}{k!} \bigl ( \varepsilon \nabla ^2 + \varepsilon ^2 \gamma (R) \bigr ) ^k
\]
are well-defined elements of $\FPDO(M,S)$.

We shall use the following two lemmas when we construct the periodic cyclic cocycle associated to the Dirac operator $\Dirac$ in Section~\ref{sec-jlo-type-cocycle}.

\begin{lemma} 
\label{lem-closure-under-ad-exp-dirac-squared}
If $T\in \PDO(U,S)$, then for every $s\in \R$  we may write 
\[ 
\operatorname{Ad}_{\exp ( s \Dirac^2 )} ( T) =\sum _{k=0}^\infty \varepsilon^k T_k
\]
where $T_k\in PDO(U,S)$ and  $\order (T_k) \le k + \order(T)$.
\end{lemma}

\begin{proof} 
This is a variation of Corollary~\ref{cor-closure-under-ad-exp-s-delta} and it is proved the same way, using the formula 
\begin{equation}
    \label{eq-bimodule-order-computation-for-dirac-squared}
\varepsilon \Delta - \Dirac ^2 = \varepsilon R_1 + \varepsilon^2 R_2 ,
\end{equation}
where $R_0$ and $R_1$ have bimodule orders $0$ and $1$, respectively. \
\end{proof}

\begin{lemma}
\label{lem-exp-dirac-squared-versus-exp-delta}
If $\Dirac$ is the Dirac operator associated to any torsion-free connection on $M$, then 
$\exp (\Dirac^2) \in \FB(M,S)$.
\end{lemma}

\begin{proof} 
The method use to prove Lemma~\ref{lem-exp-delta-prime-in-fb} applies here, in view of \eqref{eq-bimodule-order-computation-for-dirac-squared}. \end{proof}

\section{Perrot's Trace}
\label{sec-perrot-trace}
The most remarkable ingredient in Perrot's work is a supertrace functional defined on the bimodule $\FB(M,S)$  (this space  is $\Z/2\Z$-graded in the usual way, using the degree operator on differential forms).  We shall construct the supertrace in this section.  

As we did in the previous section, we shall begin by ignoring spinors: we shall first construct a trace functional on the space $\FB(M)$ (which is not $\Z/2\Z$-graded). The first ingredient in the definition of the  trace  here is the usual method of ``integrating'' homogeneous functions on $\Tcirc M$ of degree $-n$, which we shall now review.

\begin{definition}
\label{def-notation-for-densities}
If $W$ is any smooth manifold, then we shall denote by $\Density (W)$ the $C^\infty (W)$-module of smooth densities on $W$ (see \cite[pp.30-31]{BGV} or \cite[Ch.14]{Lee} for instance).
\end{definition}

View the cotangent sphere bundle $S^*M$ as the quotient of $\Tcirc M$ by the action of positive scalar multiplication on $\Tcirc M$ (generated by the Euler vector field) and denote by 
\[
p \colon \Tcirc M \longrightarrow S^*M
\]
the projection mapping.  In addition, denote by $\nu$ the volume form on $\Tcirc M$  that in any local coordinate system has the form
\[
\nu = dx^1 dx^2 \cdots dx^n d \xi_1d\xi_2 \cdots  d \xi _n.
\]

\begin{definition}
\label{def-sigma-f}
Suppose that  $f\in \Poly^{-n}(M)$.   The $(2n{-}1)$-form  $f \cdot  \iota _E\, \nu$ is basic for the action of $\R$ on $\Tcirc M$,  and we shall denote by $\sigma_f$ the  unique $(2n{-}1)$-form  on $S^*M$ such that 
\[
p^* \sigma_f =  f \cdot  \iota _E \nu  .
\]
\end{definition}

\begin{definition} 
\label{def-orientation-on-s-star-m}
Throughout the paper we shall equip $S^*M$ with the  unique orientation  such that $\int_{S^*M}  \sigma _f \ge  0$ (with $\sigma_f$ as in Definition~\ref{def-sigma-f}) for every nonnegative $f\in \Poly^{-n}(M)$. 
\end{definition} 

We obtain from  the definitions above   an isomorphism of $C^\infty (S^*M)$-modules 
\begin{equation}
    \label{eq-order-minus-n-functions-as-densities}
\Poly ^{-n}(M) \stackrel \cong \longrightarrow \Density (S^*M) .
\end{equation}
We shall use this identification throughout the rest of the paper, and for $f\in \Poly^{-n}(M)$ we shall write 
\begin{equation}
\label{eq-f-integral-definition0}
\fint _{S^* M} f =\text{integral of the associated density} =  \int_{S^*M} \sigma _{f} .
\end{equation}
For general $f\in \Poly (M)$ we shall write 
\begin{equation}
\label{eq-f-integral-definition}
\fint _{S^* M} f  = \fint_{S^*M} f^{[-n]},
\end{equation}
where $f^{[-n]}$ denotes the degree $-n$ component of $f$.

\begin{remark} 
The above integral is a commutative version of the   noncommutative residue on pseudodifferential symbols.  Indeed, up to an overall constant, the noncommutative residue of an order $-n$ pseudodifferential symbol is the integral of its principal symbol function, as defined above.
\end{remark}

The field $F$ of formal Laurent series, the algebra $\FPDO(U)$ and its subspace $FB(U)$  all carry the usual adic topologies,   in which a base for the open neighborhoods of $0$ is the collection of   spaces 
\[
\Bigl \{\, \sum _{k=\ell}^{\infty} \varepsilon ^k D_k \Bigr \} 
\]
for $\ell = 0,1,2,\dots$ (with $D_k$ either a scalar, in the case of $F$, or a partial differential operator in the other cases).  

Now denote by $\FB_{\mathrm{c}}(U)$ the subspace of $\FB (U)$ comprised of series whose coefficents are partial differential operaotrs on $\Tcirc U$ that are compactly supported in the $U$-direction, and define $\Poly_c (U)$ similarly. The main result of this section will be  as follows: 

\begin{theorem}
\label{thm-uniqueness-of-the-trace}
Let $U\subseteq M$ be an open set and with a fixed system of local coordinates.  
There is a unique   continuous \textup{(}for the adic topologies\textup{)} $F$-linear functional
\[
\Tr_F \colon \FB_{\mathrm{c}}(U) \longrightarrow  F
\]
such that 
\[
\Tr_F \bigl ( f \cdot \exp (\varepsilon \Delta ) \cdot g \bigr ) =  \varepsilon ^{-n} \fint_{S^*U} f g
\]
for all  $f, g\in \Poly_c  (U)$.  This functional has the following trace property: 
\[
\Tr_F (D X ) = \Tr_F (XD) \qquad \forall X \in \FB_c (U),\,\, \forall D \in \PDO (U).
\]
\end{theorem}

The uniquess part of the theorem is easy to prove.

\begin{proof}[Proof of uniqueness in Theorem~\ref{thm-uniqueness-of-the-trace}]
Since the \emph{finite} sums $\sum  \varepsilon ^k D_k \exp (\varepsilon \Delta)$ are dense in $\FB_c(U)$, it suffices to show that $\Tr_F$  is determined on these by the formula in the statement of the theorem.  We shall show that each such finite sum  is a finite Laurent polynomial in $\varepsilon$ with coefficients of the form $f \,  \exp (\varepsilon \Delta ) \, g$.  This certainly suffices. To do so, we need only note that  the formulas 
\begin{multline*}
 f \, \exp (\varepsilon \Delta) \,  x^j g  - f x^j \,\exp(\varepsilon \Delta ) \, g 
\\
\begin{aligned}
	& =   \varepsilon f\, \frac{\partial}{\partial \xi_j} \exp(\varepsilon \Delta ) \,  g \\
	& =  \varepsilon \frac{\partial}{\partial \xi_j}  f\cdot\exp(\varepsilon \Delta ) \, g 
-   \varepsilon \frac{\partial f}{\partial \xi_j} \exp(\varepsilon \Delta ) \, g 
\end{aligned}
\end{multline*}
and 
\[
\exp (\varepsilon \Delta) \, \xi^ i  - \xi^ i\,\exp(\varepsilon \Delta ) =   \varepsilon \frac{\partial}{\partial x^ i} \,  f\,\exp(\varepsilon \Delta ) \, g 
-   \varepsilon \frac{\partial f}{\partial x^i} \exp(\varepsilon \Delta ) \, g  ,
\]
show that the span of all $f\, \exp (\varepsilon \Delta )  \, g $ over the Laurent polynomials in $\varepsilon$ is closed under left multiplication by operators in $\PDO(U)$.
\end{proof}

Although the proof of uniqueness suggests an approach to the construction of the trace,   it is more convenient to proceed in a different direction.
Recall that if   $A$  is  a positive-definite $2n{\times}2n$ matrix, and if $p\colon \R^{2n} \to \mathbb{C}$ is a polynomial function, then
\begin{multline*}
\int _{\R^{2n}} p(w) \exp \bigl( - \tfrac 12 \langle w,  A w \rangle\bigr ) \, dw 
\\
= \frac{(2\pi )^n} {\sqrt{\det (A)}} p(\sqrt{-1}\partial_w)  \exp \bigl(-\tfrac 12 \langle w , A^{-1} w \rangle \bigr) \Bigr | _{w=0} 
\end{multline*}
(see for instance \cite[Sec.7.6]{HormanderVol1}). 
The following definition in effect extends this formula  to the matrix 
$A= \left [ \begin{smallmatrix} 0 & I \\ I & 0\end{smallmatrix} \right ]$, which is not positive-definite.

\begin{definition}
\label{def-gaussian-integral}
Let  $p\colon \R^{2n} \to \mathbb{C}$ be a polynomial function
We define  
\begin{multline*}
\Bigl \langle  p \bigl (\partial/\partial {x},\partial/\partial {\xi}\bigr ) \exp \bigl ( \varepsilon \Delta \bigr ) \Bigr \rangle
\\
\begin{aligned} 
&= 
\varepsilon ^{-n}   p\bigl (\sqrt{-1} \partial/\partial x,\sqrt{-1} \partial /\partial \xi\bigr )   \exp \bigl(\varepsilon ^{-1} x^i\xi_i \bigr ) \,\, \Bigr | _{x=0, \xi=0}  \\
& = 
\varepsilon ^{-n}   p(  \partial/\partial x, \partial /\partial \xi)   \exp \Bigl(-\varepsilon^{-1} x^i\xi_i \Bigr ) \,\, \Bigr | _{x=0, \xi=0} 
\end{aligned}
\end{multline*}
This value lies in $\C[\varepsilon^{-1},\varepsilon]$.
\end{definition}

\begin{remark} 
A comparison of the two integral formulas shows that we have dropped an overall multiplicative constant of $(2 \pi i )^n$.  For our purposes it is a little tidier to do so.
\end{remark}

Perrot calls the construction in the definition above the \emph{contraction} of $p(\partial/\partial x,\partial /\partial \xi)\exp(\varepsilon \Delta)$.  It will be convenient to follow Perrot and extend this concept of contraction  to cover more cases, as follows:

\begin{definition} 
\label{def-contraction}
If $D$ is any operator in $\PDO(U)$, say
\[
D =   \sum_{\alpha, \beta} f_{\alpha\beta }   \frac{\partial^\alpha }{ \partial{x}^\alpha}\frac{\partial^\beta}{\partial \xi_\beta },  
\]
then we define $\bigl \langle D \exp(\varepsilon\Delta ) \bigr\rangle \in \Poly(U)[\varepsilon ^{-1} , \varepsilon]$ by
\[
\begin{aligned}
\bigl \langle D  \exp \bigl ( \varepsilon \Delta \bigr ) \bigr \rangle
& = 
 \sum _{\alpha, \beta} f_{\alpha\beta} \cdot \Bigl \langle  \frac{\partial^\alpha}{\partial x ^\alpha}\frac{ \partial^\beta}{\partial \xi_\beta }  \exp \bigl(\varepsilon \Delta  \bigr )  \Bigr \rangle  \\
&  =\varepsilon ^{-n}  \sum _{\alpha, \beta} f_{\alpha\beta} \cdot \Bigl ( \frac{\partial^\alpha}{\partial x ^\alpha}\frac{ \partial^\beta}{\partial \xi_\beta }   \exp \bigl(-\varepsilon ^{-1} x^i\xi_i \bigr ) \Bigr ) \Bigr | _{x=0, \xi=0} 
 .
 \end{aligned}
\]
\end{definition}

Later on we shall make use of the following alternative formula for the trace. Form the function $\exp (\varepsilon ^{-1} (x^i-y^i)(\xi_i-\eta_i))$ on the direct product of two copies of $\Tcirc U$, and then for an operator $D \in \PDO(U)$ form the function
\[
 D \exp (-\varepsilon ^{-1} (x^i-y^i)(\xi_i-\eta_i)) 
\]
on the product, with $D$ acting on the first copy (the $x$- and $\xi$-variables) alone.  If the restriction to the diagonal of this function  is regarded as a function on $\Tcirc U$ via the projection onto the first factor of the product, then
\begin{equation}
    \label{eq-reconciliation-with-perrot-contraction}
\bigl \langle  D \exp (\varepsilon \Delta) \bigr \rangle  = 
\varepsilon^{-n} \bigl ( D \exp (-\varepsilon ^{-1} (x^i-y^i)(\xi_i-\eta_i)) \bigr ) \Big \vert _{x^i=y^i,\,\, \xi_i=\eta_i} .
\end{equation}
 
 We are now ready to give the formula for Perrot's trace functional.

\begin{definition}
\label{def-perrot-trace-scalar-version}
Let $X\in \FB (U)$.  Write $X$ as a Laurent series in $\varepsilon$
with finite singular part, 
\[
X = \sum  \varepsilon ^k D_k \exp (\varepsilon \Delta ) ,
\]
as in Definition~\ref{def-definition-of-fbu}. Assuming that each $D_k$ has compact support in the $U$-direction, we define  $\Tr_F (X) $ in the field of complex Laurent series with finite singular parts by
\begin{equation*}
\Tr_F  (X  ) 
=
\sum_{ k} \varepsilon ^ k \fint_{S^*M}  \bigl \langle D_k  \exp (\varepsilon \Delta) \bigr \rangle .
\end{equation*}

\end{definition}

Since an infinite sum is involved, we need to check that the definition actually makes sense: 

\begin{lemma}
The infinite sum  in  Definition~\textup{\ref{def-perrot-trace-scalar-version}} converges in the adic topology  on the field $F$ of formal Laurent series with finite singular parts. 
\end{lemma}

\begin{proof} 
The bimodule order plays an important role here: we shall use the fact that   if $X \in \FB(U)$, and if 
\[
X = \sum_{\alpha, \beta, p , k} \varepsilon ^ k    f_{\alpha\beta p k}   
  \frac{\partial^\alpha }{ \partial{x}^\alpha}\frac{\partial^\beta }{\partial \xi_\beta } \exp  ( \varepsilon \Delta ),
\]
then according  to   Definition~\ref{def-definition-of-fbu}, there is some $N$ for which 
\[
p  + 3 |\alpha | - | \beta | \le k + N
\]
for all $\alpha$, $\beta $, $p$ and $k$ for which $f_{\alpha\beta p k} \ne 0$.

Recall that  
\[
 \fint _{S^*M}  f_{\alpha\beta p k} \ne  0 \quad 
 \Rightarrow \quad p = -n .
 \]
In addition, it follows from Definition~\ref{def-gaussian-integral} that 
\[
\Bigl \langle 
  \frac{\partial^\alpha }{ \partial{x}^\alpha}\frac{\partial^\beta}{\partial \xi_\beta }  \exp  ( \varepsilon \Delta ) \Bigr \rangle \ne  0
  \quad \Rightarrow \quad 
\alpha \ne \beta  ,
\]
while if $\alpha = \beta$, then 
\[
\Bigl \langle 
  \frac{\partial^\alpha }{ \partial{x}^\alpha}\frac{\partial^\beta}{\partial \xi_\beta } \exp  ( \varepsilon \Delta )
  \Bigr \rangle= \text{scalar} \cdot \varepsilon ^{-|\alpha| - n}.
\]
Of course if $\alpha = \beta$, then  
\[
p  + 3 |\alpha | - | \beta | = p + 2 | \alpha | .
\]
and so if in addition $p=-n$, then 
\[
-n + 2 | \alpha | \le k + N
\]
and therefore 
\[
|\alpha| \le \frac k 2 + m ,
\]
where $m=(N+n)/2$.
We find that the series giving the trace in Definition~\ref{def-perrot-trace-scalar-version} has the form
\[
\Tr _F (X) = \sum _k \sum _{
        |\alpha| \le  k / 2 + m
}  c_{k,\alpha}\, \varepsilon^{k -|\alpha| -n} .
\]
Since $k-|\alpha|-n\ge k/2 + m - n$ under the  indicated condition on $\alpha$, the sum is certainly convergent in the adic topology, as required: the powers of $\varepsilon$ are bounded from below each power occurs only finitely many times in the sum.
\end{proof}

We shall now prove that $\Tr_F$ has the trace property in Theorem~\ref{thm-uniqueness-of-the-trace}. To do so, we shall follow the argument of Perrot closely.  The proof obviously reduces immediately to consideration of the three cases 
\begin{equation}
    \label{eq-three-cases-for-d}
D = \frac{\partial}{\partial x^i}, \quad 
D = \frac{\partial}{\partial \xi _j}\quad \text{and} \quad 
D = f .
\end{equation}
Moreover using the continuity of the trace and its $F$-linearity, the proof simultaneously reduces to the case where 
\[
X = g\, \frac{\partial^\alpha}{\partial x^\alpha} \frac{\partial^\beta }{\partial \xi_\beta} \exp (\varepsilon \Delta).
\]
In the first two cases from \eqref{eq-three-cases-for-d} we have 
\[
[D,X] = \frac{\partial g}{\partial x^i} \frac{\partial^\alpha}{\partial x^\alpha} \frac{\partial^\beta }{\partial \xi_\beta} \exp (\varepsilon \Delta)
\quad\text{and}\quad 
[D,X] = \frac{\partial g}{\partial \xi_j} \frac{\partial^\alpha}{\partial x^\alpha} \frac{\partial^\beta }{\partial \xi_\beta} \exp (\varepsilon \Delta) .
\]
To handle these two cases of the proof it suffices to show that 
\[
\fint _{S^* U}  \frac{\partial g}{\partial x^i} = 0 \quad \text{and} \quad  \fint _{S^* U} \frac{\partial g}{\partial \xi_j} = 0 .
\]
In the first case we may as well assume that $g$ has degree $-n$ (otherwise the integral is certainly zero, by definition), and then by Stokes' theorem
\[
\fint _{S^* U}  \frac{\partial g}{\partial x^i} = 
\int _{S^* U} d\, \iota _{\partial/\partial x_i} \sigma _g = 0 ,
\]
in the notation of Definition~\ref{def-sigma-f}.  In the second case we may as well assume that $g$ has degree $1{-}n$.  But then the form $\iota_E \iota _{\partial/\partial \xi _j} g \nu$  on $\Tcirc U$ is basic, so it is the pullback of a form $\tau_g$ on $S^* U$, and  from 
\[
d ( \iota_E \iota _{\partial/\partial \xi _j} g \nu) = \iota _E \frac{\partial g }{\partial \xi_j} \nu,
\]
we find that 
\[
\fint _{S^* U}  \frac{\partial g}{\partial \xi_j} = 
\int _{S^* U}   \sigma _{\partial g/\partial \xi_j} = \int _{S^* U} d \tau_g = 0 ,
\]
again by Stokes' theorem.

It remains to consider the case $D = f$ in \eqref{eq-three-cases-for-d}. Consider first the cases where $f$ is a coordinate function: 
$f = x^i$ or $f = \xi _j$.

\begin{lemma} 
\label{lem-explicit-commutator-computation-for-trace-property}
For any $i$, and $j$ and any $\alpha$ and $\beta$, 
\[
\Bigl \langle \bigl [x^i, \frac{\partial^\alpha}{\partial x^\alpha}
\frac{\partial^\beta}{\partial \xi_\beta } \exp (\varepsilon \Delta )\bigr ] \Bigr \rangle = 0
\quad \text{and} \quad 
\Bigl \langle \bigl [\xi_j, \frac{\partial^\alpha}{\partial x^\alpha}
\frac{\partial^\beta}{\partial \xi_\beta} \exp (\varepsilon \Delta )\bigr ] \Bigr \rangle = 0.
\]
\end{lemma}

\begin{proof}
We shall present the argument for $x^i$ (the argument for $\xi_j$ is exactly the same) and we'll describe  the case of $n{=}1$ dimensions (higher dimensions present additional notational complexities but are otherwise identical).

We calculate that 
\[
 \bigl [ \frac{\partial^a}{\partial x^a}
\frac{\partial^b}{\partial \xi ^b} \exp (\varepsilon \Delta ),x \bigr ] = 
a \frac{\partial^{a-1}}{\partial x^{a-1} }\frac{\partial^b}{\partial \xi ^b} \exp (\varepsilon \Delta ) +\varepsilon  \frac{\partial^a}{\partial x^a}
\frac{\partial^{b+1}}{\partial \xi ^{b+1}} \exp (\varepsilon \Delta )
\]
and therefore that the contraction is 
\[
\frac{\partial^b}{\partial \xi ^b}\Bigl ( a \frac{\partial^{a-1}}{\partial x^{a-1} } + \varepsilon \frac{\partial^{a}}{\partial x^{a} }\frac{\partial}{\partial \xi } \Bigr )  \exp ( - \varepsilon^{-1} x \xi ) \Big \vert _{x=0,\xi=0}.
\]
But 
\[
\Bigl ( a \frac{\partial^{a-1}}{\partial x^{a-1} } + \varepsilon \frac{\partial^{a}}{\partial x^{a} }\frac{\partial}{\partial \xi } \Bigr )  \exp ( - \varepsilon^{-1} x \xi ) = (-1)^{a-1} \varepsilon ^{-a}\xi^{a-1} x \exp(-\varepsilon^{-1} x \xi) .
\]
So, thanks to the factor of $x$, the contraction is $0$.
\end{proof}

\begin{remark} 
The lemma can be understood more conceptually (but at the expense of a longer argument) using the fact that the contraction is derived from an integral formula on $\R^{2n}$, and the fact that the integral of a partial  derivative in that context is zero.
\end{remark} 

\begin{proof}[Proof of the trace property in Theorem~\ref{thm-uniqueness-of-the-trace}]
The calculations in the proof of Lemma~\ref{lem-explicit-commutator-computation-for-trace-property} show that if $p  $ is any polynomial on $\R^{2n}$, possibly with coefficients that are  polynomials in $\varepsilon$, then 
\[
\Bigl [x^i, p \Bigl (\frac{\partial }{\partial x },
\frac{\partial }{\partial \xi}\Bigr ) \exp (\varepsilon \Delta )\Bigr ]
=q\Bigl (\frac{\partial}{\partial x} ,
\frac{\partial }{\partial \xi }\Bigr )  \exp (\varepsilon \Delta )
\]
for some other $q$.  The same is true for commutators with $\xi_j$, and as a result the lemma in fact shows that if $\mu$ and $\nu$ are any multi-indices, then 
\[
\Bigl \langle 
\ad_x^\mu \ad_\xi^\nu \Bigl ( p \Bigl (\frac{\partial }{\partial x },
\frac{\partial }{\partial \xi}\Bigr ) \exp (\varepsilon \Delta )\Bigr)
\Bigr\rangle = 0.
\]
But if $f$ is any polyhomogeneous function, then 
\begin{multline*}
\Bigl [  p \Bigl (\frac{\partial }{\partial x },
\frac{\partial }{\partial \xi}\Bigr ) \exp (\varepsilon \Delta ), f\Bigr ]
\\
=
\sum_{\mu,\nu\ge 0}
\frac{1} {\mu! \nu!} \frac{\partial^\mu\partial^\nu f  }{\partial x^\mu\partial \xi^\nu} \cdot \ad_x^\mu \ad_\xi^\nu \Bigl ( p \Bigl (\frac{\partial }{\partial x },
\frac{\partial }{\partial \xi}\Bigr ) \exp (\varepsilon \Delta )\Bigr),
\end{multline*}
with the sum convergent in the adic topology.  So the contraction of the commutator is $0$. We have handled the third and final case in \eqref{eq-three-cases-for-d}, and therefore the proof is complete.
\end{proof} 

The second crucial fact about $\Tr_F$ is its coordinate independence: 

\begin{theorem}
\label{thm-independence-of-the-trace}
The trace functional  $\Tr_F$ is independent of the choice of coordinate system used to define it.
\end{theorem}

\begin{proof}
The theorem isn't an immediate consequence of the uniqueness part of Theorem~\ref{thm-uniqueness-of-the-trace} because the Laplacian $\Delta$ is not invariant under coordinate changes.  What we need to show is that if $\Tr_F$ and $\Delta$ are constructed using a first coordinate system, and if $\Delta'$ is the Laplacian in a second coordinate system, then
\[
\Tr_F ( f \exp (\varepsilon \Delta)) =\fint_{S^*M} f = \Tr_F (f \exp( \varepsilon \Delta ')).
\]

As we already noted in the proof of Lemma~\ref{lem-exp-delta-prime-in-fb}, the second Laplacian has the form 
\[
\Delta ' = \Delta  + a_i \frac{\partial}{\partial \xi_i} + b^{j}_{k\ell} \xi_j  \frac{\partial^2}{\partial \xi_k \partial \xi_\ell}
\]
in the first coordinate system, where $a_i$ and $b^j_{k\ell}$ are smooth functions on $M$.  If $T = \Delta ' -\Delta $ then we can write (using  the contraction operation  associated to the first coordinate system)
\[
\bigl \langle f \exp(\varepsilon \Delta ') \bigr \rangle  - 
\bigl \langle f \exp (\varepsilon \Delta) \bigr \rangle
= \int ^1_0  \bigl \langle f \exp (s\varepsilon \Delta) \varepsilon T \exp ((1{-}s)\varepsilon \Delta ')\bigr \rangle \, ds .
\]
If we expand the operators $\exp (s\varepsilon \Delta) \varepsilon T \exp ((1{-}s)\varepsilon \Delta ')$ as power series in $\varepsilon$, then the coefficient of $\varepsilon^k$, which is an element of $\PDO (U)$, has differential order $k{-}1$ in the $x$-direction (because $T$ includes no $x$-derivative) and is homogeneous of degree $-k$ in the vertical direction (since each of $\Delta$, $\Delta'$ and $T$ has degree $-1$).  This means that if we write the same operators in the standard form 
\[
\exp (s\varepsilon \Delta) \varepsilon T \exp ((1{-}s)\varepsilon \Delta ') = 
\sum _{k=0}^\infty \varepsilon^k D_k \exp (\varepsilon \Delta)
\]
and then write 
\[
D_k = \sum g_{k,\alpha,\beta} \frac{\partial^\alpha}{\partial x^\alpha}  \frac{\partial^\beta}{\partial \xi_\beta} ,
\]
then 
\[
|\alpha| \ge  |\beta| \quad \Rightarrow \quad g_{k,\alpha,\beta} = 0.
\]
The contraction is therefore zero and the theorem is proved.
\end{proof}

With Theorems\ref{thm-uniqueness-of-the-trace} and \ref{thm-independence-of-the-trace} available, we can now extend the trace beyond coordinate charts to a continuous and $F$-linear trace functional
\[
\Tr_F \colon \FB (M) \longrightarrow F 
\]
in the obvious way using partitions of unity.
But for the purposes of this paper we only need a part of the trace we have just defined, which is scalar-valued, not $F$-valued: 

\begin{definition}
\label{def-perrot-trace}
Let $M$ be a smooth, closed manifold. The \emph{Perrot trace} on $\FB(M)$ is 
\begin{gather*}
    \Tr \colon \FB (M) \longrightarrow \C \\
    \Tr (X) = \text{Coefficient of $\varepsilon ^0$ in $\Tr_F(X)$}.
\end{gather*}
\end{definition}

Finally, we need to pass from a (scalar)  trace on $\FB(M)$ to a supertrace on $\FB(M,S)$.  First denote by $\Gamma  (\operatorname{End}(S) ) $ the  smooth polyhomogeneous sections of the endomorphism bundle of $S$.  This is a $\Z/2\Z$-graded algebra over $\Poly (M)$ and there is a standard $\Poly (M)$-linear supertrace 
\[
\str_\Lambda \colon \Gamma (\operatorname{End}(S)  ) 
\longrightarrow \Poly(M)
\]
In local coordinates: 
\[
\str_\Lambda (\psi^\mu \overline \psi _\nu ) 
= \begin{cases}  1 & |\mu| = |\nu| =n\\ 0 & \text{otherwise}.
\end{cases} 
\]
Then, in any coordinate chart $U\subseteq M$ we  may use the  natural (coordinate-dependent) vector space  isomorphism 
\[
   \FB(U,S) \cong \Gamma  (\operatorname{End}(S\vert _U))  \otimes _{\Poly (U)} \FB (U)  
\]
to define 
\begin{gather*}
\STr_F \colon \FB(U,S) \longrightarrow  F \\
\STr_F (A \otimes X ) = \Tr ( \str_\Lambda (A)X)
\end{gather*}
(strictly speaking this is defined only on elements compactly supported in the $U$-direction).  Thanks to Theorem~\ref{thm-uniqueness-of-the-trace}, the above is a bimodule supertrace for the left and right actions of the algebra 
\[
\PDO (U,S)\cong \Gamma  (\operatorname{End}(S\vert _U))  \otimes _{\Poly (U)} \PDO (U).
\]
Thanks to Theorems~\ref{thm-uniqueness-of-the-trace} and \ref{thm-independence-of-the-trace} it is independent of the choice of coordinates used to define it.

We can now use partitions of unity to define (independently of all choices)  a supertrace
\[
\STr_F \colon \FB (M) \longrightarrow F ,
\]
and finally we make the following definition: 

\begin{definition}
\label{def-perrot-supertrace}
Let $M$ be a smooth, closed manifold. The \emph{Perrot supertrace} on $\FB(M,S)$ is 
\begin{gather*}
    \STr \colon \FB (M,S) \longrightarrow \C \\
    \STr (X) = \text{Coefficient of $\varepsilon ^0$ in $\STr_F(X)$}.
\end{gather*}
\end{definition}

\section{The JLO-Type Cocycle}
\label{sec-jlo-type-cocycle}

The purpose of this section is to construct from the Dirac operator $\Dirac$, together with Perrot's supertrace  and one more ingredient, a periodic cyclic cocycle for the algebra $C^\infty (S^*M)$. 

Smooth functions on $S^*M$ can be viewed as order-one homogeneous functions on $\Tcirc M$. Seen this way, the algebra $C^\infty (S^*M)$, together with bimodule $FB(M,S)$, equipped with its supertrace, and the operator $\Dirac$ form something similar to one of Connes'  \emph{spectral triples}   \cite[Sec.IV.2]{ConnesNCGBook}\footnote{In this reference the older term \emph{$K$-cycle} is used, rather than \emph{spectral triple}.}:    the bimodule, with its trace, plays a role similar to the role of the Hilbert space in a spectral triple, with its associated operator trace (or perhaps with a residue trace obtained from this operator trace; see for instance \cite{Higson04} for an exposition).

With this in mind we shall   attempt to adapt the well-known  JLO-cocycle in periodic cyclic theory to Perrot's context. 
Let $\Sigma^p$ be the standard $p$-simplex. If  $s\in \Sigma^p$ and  $X^0,\dots , X^p\in \PDO (M,S)$, then let us write 
\begin{equation}
\label{eq-def-of-angle-bracket-jlo-term}
\bigl \langle 
X^0,\dots, X^{p}
\bigr \rangle_{s} 
 =
X^0\exp(s_0 \Dirac^2)X^1\exp(s_1 \Dirac^2)\cdots X^{p} 
\exp(s_{p} \Dirac^2) .
\end{equation}
The most obvious  JLO cocycle that one might try to attach to $\Dirac$ is given by the formula 
\begin{equation}
\label{eq-wrong-jlo-cocycle}
JLO^{\Dirac}(a^0,\dots, a^p) =
\int_{\Delta ^{p}}
\STr \Bigl (
\bigl \langle a^0, [\Dirac, a^1], \dots   , [\Dirac, a^p ] \bigr \rangle_s 
\Bigr ) \, ds .
\end{equation}
This is the original formula given by Jaffe, Lesniewski and Osterwalder \cite{JLO88}; see \cite[Sec.8]{Quillen88} or \cite[Sec.IV.8]{ConnesNCGBook} for expositions.

However it follows from the supertrace property that the  cocycle \eqref{eq-wrong-jlo-cocycle} vanishes in odd degrees, whereas the Todd class that we are trying to understand corresponds to a cyclic cocycle that vanishes in even degrees.  So the above is not the cocycle we are seeking.\footnote{In fact it may be shown that \eqref{eq-wrong-jlo-cocycle} vanishes in \emph{all} degrees. In the more elaborate context of Perrot's work on pseudodifferential symbols, this corresponds to the vanishing of the noncommutative residue on pseudodifferential projectors \cite[Thm.6.5]{Perrot2013}.}
 
To address this issue we need one further ingredient,   as follows.
Fix a Riemannian  metric on $M$, and hence a Euclidean structure on the bundle $T^*M$.  Define 
\[
q\colon \Tcirc M \longrightarrow (0,\infty),\qquad q(\alpha) = \| \alpha\|,
\]
and then define a derivation
\[
\delta \colon \PDO(M,S) \longrightarrow \PDO(M,S)
\]
by the formula 
\[
\delta (X) = \ad_{\log (q)} (X)
\]
(as usual we are adapting Perrot's work here, although we are adapting it to our simplified context).  Note that 
\[
\delta \Bigl (\frac{\partial}{\partial x^i }\Bigr )  = - \frac{\partial q } { \partial x^i} \cdot \frac{1}{q}
\quad 
\text{and} 
\quad 
\delta \Bigl (\frac{\partial}{\partial \xi_j }\Bigr )  = - \frac{\partial q } { \partial \xi_j} \cdot \frac{1}{q},
\]
while $\delta$ vanishes on other (local) generators of $\PDO(M,S)$.  So $\delta$ does indeed map $\PDO(M,S)$ into itself, despite the fact that $\log (q)$ is not a polyhomogeneous function.

Now $\delta$ extends to a derivation of the algebra $\FPDO(M,S)$ by applying $\delta$ to each coefficient in a formal series.  And this extension maps the bimodule  $\FB(M,S)$ to itself.  This follows from the formula 
\[
\begin{aligned}
\delta ( \exp (\varepsilon \Delta)) 
    & = \int _0^1 \exp (s \varepsilon \Delta) \varepsilon \delta (\Delta)\exp((1{-}s)\varepsilon \Delta) \, ds 
\\
    & = \int _0^1 \operatorname{Ad}_{\exp (s \varepsilon \Delta)}\bigl (\varepsilon \delta (\Delta) \cdot \exp(\varepsilon \Delta) \bigr )   \, ds 
\end{aligned}
\]
together with Lemma~\ref{cor-closure-under-ad-exp-delta}.

\begin{lemma}
\label{lem-closed-trace}
If $X\in \FB(M,S)$, then  
$\STr (\delta (X)) = 0$. 
\end{lemma}

\begin{proof} 
Let $X$ be any linear partial differential operator on $\Tcirc M$ and let $u\in \R$.  Consider the binomial-type formula 
\[
\ad_{q^u} (X) = \sum _{k=1}^\infty \binom{u}{k} \ad_q^k(X) q^{u-k} ,
\]
which actually involves a finite sum on the right-hand side, since if $h$ is any smooth function on $\Tcirc M$, then the action of $\ad_h$ is locally nilpotent on partial differential operators.
Differentiating the formula with respect to $u$ and then setting $u=0$, we obtain 
\[
\ad_{\log (q)}( X) = \sum _{k=1}^\infty \frac {1}{k} (-1)^{k-1} \ad_q^k(X) q^{-k}
\]
Again the right-hand side is a finite sum, and it is a combination of commutators, since $\ad_q^k(X)q^{-k}=\ad_q^k(Xq^{-k})$ for all $k$.

Now we turn to $\FPDO(M,S)$. The above formulas involve infinite sums when $X\in \FPDO (M,S)$, but they are convergent in the adic topology and for that reason remain valid. Restricting to $\FB(M,S)$, it follows that every element of the form
$\delta (X)=\ad_{\log(q)}(X)$ is a limit of sums of commutators, and since the trace is continuous the trace of $\delta (X)$ is zero.
\end{proof}

\begin{remark}
It is worth noting that the trace does \emph{not} extend as a trace if we adjoin $\log (q)$ to the bimodule $\FB(M)$ (so this is not the explanation for the lemma above). For instance, in case where $M$ is a circle, if $q=|\xi|$, and if $s$ is the characteristic function of $\{\,\xi > 0\,\}$ (which is a polyhomogeneous function), then 
\[
\Bigl [ \frac {\partial}{\partial \xi}, s \log (q) \exp (\varepsilon \Delta)\Bigr]
= s \xi^{-1} \exp (\varepsilon \Delta) ,
\]
and the trace of the right-hand side is nonzero (despite it being a commutator in the enlarged bimodule).
\end{remark}

\begin{remark}
On a related note, although the definition of $\delta$ requires a choice of Riemannian metric, different choices yield derivations that differ only by an inner derivation.
\end{remark}

Using Lemma~\ref{lem-closed-trace} it is not difficult to adjust the definition of the JLO cocycle so as to incorporate $\delta$.  The following theorem is a generalization of the well-known fact that if $\delta$ is any derivation on an algebra, and if $\tau$ is a trace with $\tau(\delta (a))=0$ for all $a$, then the formula $\varphi (a^0,a^1) = \tau (a^0\delta (a^1))$ defines a cyclic cocycle.  As Rodsphon explains in \cite{Rodsphon2015}, the theorem fits very naturally into Quillen's formalism for the JLO cocycle \cite{Quillen88}, and this is probably the best way of understanding it.

First a simple preliminary calculation: 

\begin{lemma} 
Let $\nabla$ be a torsion-free connection on $M$ and let $\Dirac$ be the associated Dirac operator, as in Definition~\textup{\ref{def-dirac-operator}}. If $p\ge 0$,  $s\in \Sigma^{p+1}$ and $X_0,\dots, X_{p+1}\in \PDO (M,S)$, then
\[
\bigl \langle 
X^0,\dots, X^{p+1}
\bigr \rangle_{s} \in \FB(M,S),
\]
to use  the notation of  \textup{\eqref{eq-def-of-angle-bracket-jlo-term}}. 
\end{lemma} 

\begin{proof} 
This is a combination of Lemmas~\ref{lem-closure-under-ad-exp-dirac-squared} and \ref{lem-exp-dirac-squared-versus-exp-delta}.
\end{proof} 

\begin{theorem}
\label{thm-jlo-type-cocycle}
The formulas 
\begin{multline*}
JLO^{\Dirac,\delta}(a^0,\dots, a^p) =
\sum_{k=1}^{p+1} (-1)^k 
\int_{\Sigma ^{p+1}}
\STr \Bigl (
\bigl \langle a^0, [\Dirac, a^1], \dots \\
\dots , [\Dirac, a^{k-1}], \delta (\Dirac), [\Dirac, a^{k}], \dots , [\Dirac, a^p ] \bigr \rangle_s 
\Bigr ) \, ds
\end{multline*}
for $p$ odd and $a^0,\dots, a^p\in C^\infty (S^*M)$
define a periodic cyclic cocycle for the algebra $C^\infty(S^*M)$.
\end{theorem}

\begin{remark}
The integrands in the defining formula for $JLO^{\Dirac, \delta}$ are polynomial functions of $s\in \Sigma^{p+1}$, so the integrals are certainly well-defined, and moreover these polynomial functions are identically zero for $p$ sufficiently large (this will be made clear in the next section), so that only finitely many components of $JLO^{\Dirac, \delta}$ are nonzero, as required in the definition of periodic cyclic cocycle.  This is in contrast to the form of the usual JLO cocycle (involving the Levi-Civita Dirac operator and the operator trace), which is rather an \emph{entire} cyclic cocycle \cite[Ch.4, Sec.7]{ConnesNCGBook}.
\end{remark}

\begin{remark} 
In general the Quillen formalism would produce a formula with further terms, involving factors $\delta (a^j)$; see \cite[Sec.3]{Rodsphon2015}. But in the case at hand, $\delta(a^j)=0$.
\end{remark}

We shall not prove  Theorem~\ref{thm-jlo-type-cocycle} in this paper, firstly because  it is a general result that is encompased by the Quillen formalism and  is explained perfectly well by Rodsphon in \cite{Rodsphon2015}, and secondly because  we shall not actually need the result: we shall calculate all the components of $JLO^{\Dirac, \delta}$, and it will be evident from these computations that they constitute a periodic cyclic cocycle.

\section{Perrot Order of a  Differential Operator} 
\label{sec-perrot-order}
The remainder of the paper will be dedicated to the computation  of the JLO-type cocycle that was  defined in  the previous section.  We shall introduce a new notion of order on differential operators that we shall call the \emph{Perrot order} (it is distinct from the bimodule order discussed earlier). The Perrot order
has the property that operators of negative order have vanishing trace. Roughly speaking, upon analyzing the constituents of the cocycle, we shall find that many parts have negative order, and so they can be removed without altering the value of the cocycle.

The Perrot order resembles somewhat the notion of order that Getzler uses to compute the index of the spinor Dirac operator \cite{Getzler83} (see \cite{RoeAsymptoticMethods} or \cite{HigsonYi19} for surveys).  Both have the property that the curvature operator\footnote{In Perrot's case, we are referring here to the operator    $\gamma(R)$ that was introduced in Definition~\ref{def-gamma-of-r} and appears in the formula for $\Dirac^2$. The actual curvature operators $R(X,Y)$ have order $1$; see the proof of  Lemma~\ref{lem-covariant-derivatives-commute-to-highest-order}.}  acquires  the same order as the square of the Dirac operator, which leads to the curvature contributing to leading order computations.  But the leading order is $2$ for Getzler, while it is  $0$ for Perrot. An important consequence is that the formal exponential $\exp (\Dirac ^2)$ also has Perrot order $0$. Indeed all of the constituents of the JLO cocycle introduced in the previous section have Perrot order $0$. 

Despite its role in capturing curvature,  the Perrot order does not in fact depend on the choice of connection, whereas Getzler's certainly does (which is a  paradox, perhaps, but not a contradiction).    But to   begin with, to define the Perrot order it is convenient to fix  a torsion-free connection on $M$ (and then later we shall show that the Perrot order does not depend on this choice).

\begin{definition} 
\label{def-perrot-filtration}
The \emph{Perrot order} of an operator in $\PDO(M,S)$ (which is an integer, possibly negative) is defined by the following:

\begin{enumerate}[\rm (i)]
\item If $f$ is any degree $p$-homogeneous smooth function on $\Tcirc M$, then $\perrotorder(f)\le p$.  

\item If $X$ is any vector field on $M$, with horizontal lift $X^H$ on $\Tcirc M$,  then  $\perrotorder (\nabla_{X^H}) \le 1$.

\item If $Y$ is any vertical vector field on $\Tcirc M$ (that is, if $Y$ commutes with functions pulled back from $M$), and if $Y$ is translation-invariant on each cotangent fiber, then $\perrotorder (\nabla _Y)\le -1$.

\item If $\omega$ is any $1$-form on $M$,  acting  on sections of $S$  over $\Tcirc M$ by left exterior multiplication, then $\perrotorder (\omega) \le 0$.

\item If $X$ is any vector field on $M$, and if $\iota (X)$ denotes the action  on sections of $S$  by contraction, then $\perrotorder (\iota (X))\le 1$.

\end{enumerate}
\end{definition}

\begin{remark}
Compare Proposition/Definition~\ref{prop-def-bimodule-order-scalar-case}, which explains the way in which an increasing filtration on $\PDO (M,S)$ is constructed from  inequalities such as those above. As with the bimodule order considered there, the inequalities in the definition above are actually equalities;  this will be clear from the lemma below.
\end{remark}

An operator has Perrot order $p$ or less if and only if its restriction to every $\Tcirc U$, with $U\subseteq M$ open, has Perrot order $p$ or less.  So the Perrot order may be computed locally.  The following lemma shows that moreover the Perrot order  is easily determined in local coordinates. It also shows that the Perrot order  is independent of the choice of affine connection $\nabla$. 

\begin{lemma}
\label{lem-order-bound-by-coordinate-filtration}
Let   $x_1,\dots , x_n$ be coordinates on an open set $U \subseteq M$ and let $p\in \Z$.  If the bundle $S$ is trivialized over $\Tcirc U$ using the standard frame for the exterior algbra bundle associated to these coordinates, then the linear space of all operators of Perrot order $p$ on less on $\Tcirc U$ is precisely the  linear span of operators of the form
\[
f  \psi^\mu \overline \psi _\nu \frac{\partial^\alpha}{\partial x^\alpha}\frac{\partial^\beta}{\partial \xi_\beta}
\]
with $f$ a homogeneous smooth scalar function on $\Tcirc U$ and with
\[
\degree (f)  +|\alpha| -  | \beta | + |\nu|   \le p.
\]
Here $\alpha$ and $\beta$ are multi-indices with non-negative integer entries, while $\mu$ and $\nu$ are multi-indices with $0$-$1$ entries.
\end{lemma}

\begin{proof}
It follows from \eqref{eq-formula-for-induced-connection} that 
\[
\frac{\partial}{\partial x^i} =  \nabla _{(\partial/\partial x^i)^H}   - \Gamma_{i\ell}^k\, \xi_k\, \frac{\partial}{ \partial \xi _\ell} + \Gamma_{i\ell}^k \,\psi^\ell\, \overline \psi_k .
\]
All   three terms on the right-hand side have Perrot order $1$ or less, and therefore $\partial / \partial x^i$ has Perrot order $1$ or less.   In addition 
\[
\frac{\partial}{\partial \xi_j} = \nabla _{\partial/\partial \xi _j} .
\]
 It follows that all the operators in the statement of the lemma have Perrot order $p$ or less.

In the reverse direction, if we define 
$ \PDO_p (U,S)$ to be the linear span of the operators in the statement of the lemma, then we obtain an algebra filtration 
for which the associated order satisfies the relations in Definition~\ref{def-perrot-filtration}.    It therefore follows from the construction of the Perrot filtration and order (\emph{c.f.} Propositon/Definition~\ref{prop-def-bimodule-order-scalar-case}) that every operator of Perrot order $p$ or less belongs to $\PDO_p(U,S)$.
\end{proof}

 Let us now extend the Perrot order to the algebra $\FPDO (M,S)$ and to the bimodule $\FB(M,S)$ as follows: a formal series $\sum \varepsilon ^k T_k$ has Perrot order $p$ or less if and only if each $T_k$ has Perrot order $p$ or less.  (With this definition, there are many elements of infinite Perrot order, but that will not be an issue for us.)  

\begin{proposition}
Let $M$ be a closed manifold.
The supertrace vanishes on every element of $\FB(M,S)$  of negative Perrot order.
\end{proposition}

\begin{proof} 
It suffices to consider negative-order generating elements 
\[
f  \psi^\mu \overline \psi _\nu \frac{\partial^\alpha}{\partial x^\alpha}\frac{\partial^\beta}{\partial \xi_\beta} \exp (\varepsilon \Delta) 
\]
of the bimodule $\FB (U,S)$ in a coordinate neighborhood.  According to the definition of the supertrace and the formula for the contraction operation in Definition~\ref{def-gaussian-integral}, if $\alpha \ne \beta$, then the supertrace of this element is zero. In addition, if  $|\nu| \ne n$ then the supertrace is zero. If $\alpha= \beta$ and if $|\nu | = n$, then the Perrot order is equal to $\degree (f) - n$, and if this is negative, then by the definition of the integral $\fint f$ in  \eqref{eq-f-integral-definition} the supertrace is zero.
\end{proof}

Our purpose in the remainder of this section is to develop a technique to compute the supertrace of order zero operators. In the lemma below we shall use the fact that if $\pi\colon W\to M$ is any submersion with compact fibers, then there is an induced morphism of $C^\infty (M)$-modules  
\[
\pi_*\colon \operatorname{Dens}(W)\longrightarrow \operatorname{Dens} (M)
\]
that is characterized by the formula 
\[
\int_M  \pi_*(\alpha ) = \int_W \alpha 
\]
(in the context of oriented smooth manifolds this is the operator of integration over the fiber). 

\begin{lemma} 
There is a unique morphism of $C^\infty (M)$-modules  
\[
\str\colon \FB(M) \longrightarrow  \operatorname{Dens}(M)
\]
such that 
\[
\STr (X) = \int_{M} \str (X)
\]
for every $X\in \FB^0(M)$.
\end{lemma}

\begin{proof} 
 The trace was constructed (locally, at first) as the integral of a density on $S^*M$, and pushing this density forward along the projection $S^*M\to M$ we obtain from a partition of unity argument the existence of a map as in the statement of the lemma. 

To prove uniqueness, if $\str$ and $\str'$ are two morphisms of modules  as in the statement of the lemma, then from the identities
\begin{multline*}
\int_{M} f\cdot \str (X) = \int_{M}   \str (f\cdot X) =  \STr (X) 
\\ 
= \int_{M}   \str' (f\cdot X) = \int_{M} f\cdot \str' (X)
\end{multline*}
we conclude that 
\[
\int _M f\cdot \bigl ( \str (X) - \str'(X)\bigr ) = 0
\]
 for all $f\in C^\infty (M)$, and hence $\str(X)-\str'(X)=0$.
\end{proof}


We shall now evaluate the density $\str (X) \in \Density (M)$ defined above  at a single point $m\in M$ under the assumption that  $X\in \FB (M,S)$ has order zero.  

Form the associated graded algebra 
\[
\gr \PDO (M,S) = \bigoplus _p \PDO _p (M,S)\big / \PDO _{p-1} (M,S)
\]
for the filtration by Perrot order. Notice that since smooth functions on $M$ (viewed as operators on sections of $S$ by pointwise multiplication) have Perrot order zero, the associated graded algebra is in fact an algebra over $C^\infty (M)$.  We can therefore speak of the fiber 
\[
\gr \PDO (M,S) \big \vert _m 
\]
at a given point $m\in M$, which is an algebra in its own right.  Our first aim is to compute this fiber.  

\begin{lemma} 
\label{lem-covariant-derivatives-commute-to-highest-order}
If $X$ and $Y$ are vector fields on $M$, then the degree $1$ classes of $\nabla_{X^H}$ and $\nabla_{Y^H}$ commute with one another in the associated graded algebra.
\end{lemma}

\begin{proof}
Let $\alpha \in \Tcirc M$ and let $m=\pi(\alpha)$.
Because the connection we are using on $\pi^* \Lambda^* T^*M$   is a pullback,  the curvature operator  $R(X^H,Y^H)\vert _\alpha $  
is  the curvature operator $R(X,Y)\vert _m$ for the original connection: 
\begin{equation}
    \label{eq-commuting-two-nablas}
  \nabla _{X^H} \nabla _{Y^H} - \nabla _{Y^H} \nabla _{X^H}  = \nabla _{[X^H,Y^H]} + R(X,Y) 
\end{equation}
(this is an identity of operators acting on the sections of $S = \pi^* \Lambda ^* T^*M$). 

If we introduce local coordinates on $M$ and write the curvature operator on $TM$ as
\[
R \Bigl(\frac{\partial}{\partial x^i}, \frac{\partial}{\partial x^j}\Bigr )\, \frac{\partial}{\partial x^\ell} =   R^k_{ij\ell}\, \frac {\partial}{\partial x^k} ,
\]
as in Definition~\ref{def-gamma-of-r}, then the curvature operator as it appears in \eqref{eq-commuting-two-nablas} is 
\[
R \Bigl(\frac{\partial}{\partial x^i}, \frac{\partial}{\partial x^j}\Bigr ) = \psi^\ell \overline \psi_k R^k_{ij\ell} ,
\]
and in particular the operator $R(X,Y)$ in \eqref{eq-commuting-two-nablas} has Perrot order $1$.

In addition,   according to   Lemma~\ref{lem-horiz-vector-fields-and-curvature},
\[
[X^H,Y^H] = [X,Y]^H + \gamma \bigl (R(X,Y)\bigr ) .
\]
The covariant derivatives attached to the two terms on the right have Perrot-orders $1$ and $0$, respectively.  So we find  from \eqref{eq-commuting-two-nablas} that 
\begin{multline*}
\perrotorder \bigl ( \nabla _{X^H} \nabla _{Y^H} - \nabla _{Y^H} \nabla _{X^H}\bigr ) 
\\
\begin{aligned}
& \le 1 \\
&
< \perrotorder \bigl ( \nabla _{X^H} \bigr ) + \perrotorder \bigl ( \nabla _{Y^H} \bigr ),
\end{aligned}
\end{multline*}
which proves the lemma.
\end{proof}

\begin{remark}
The lemma stands in contrast to the situation with the Getzler's order, in which the commutator covariant derivatives is a curvature operator.
\end{remark}

Other, more easily derived,  relations in the associated graded algebra are as follows: 
\begin{enumerate}[\rm (i)]

\item  The degree $-1$ classes of the operators $\nabla_Y$ associated to vertical vector fields that are fiberwise translation invariant (that is, scalar combinations of the $\nabla_{\partial /\partial \xi_j}$) commute with one another and the classes of the $\nabla_{X^H}$.

\item The degree $0$ and $1$ (respectively) classes of the exterior multiplication and contraction operators anticommute with one another and commute with all the classes of covariant derivatives already mentioned. 

\item The degree $p$ classes of scalar degree $p$ functions commute with one another and with all the classes mentioned above, except that \[
[\nabla_{\partial/\partial \xi _j}] [f] - [f] [\nabla_{\partial / \partial \xi _j}] = [\partial f / \partial \xi _j]
\]
as in $\PDO(M,S)$.
\end{enumerate}

\begin{definition} 
\label{def-symbol-algebra-at-m}
For $m\in M$, denote by $\SymbolAlg  _m$   the tensor product of the symmetric algebra of $T_mM$ (whose elements we shall write as constant coefficient differential operators on $T_mM$), the algebra of polyhomogeneous partial differential operators on $T^*_mM$ and the exterior algebra $\Lambda ^* (T_mM \oplus T^*_m M)$ (which, given local coordinates on $M$ near $m$, we shall regard as generated by $dx^i\in T^*_mM$ and $d \xi_j\in T_mM$). 
\end{definition} 

We shall suppress tensor product signs when describing elements of $\SymbolAlg_m$, such as for example
\begin{equation}
\label{eq-example-of-symbol-operator}
dx^i \frac{\partial}{\partial x^i} + d \xi_j \frac{\partial }{\partial \xi _j}
\in \SymbolAlg_m
\end{equation}
We equip $\SymbolAlg_m$ with an integer grading  by assigning degree $1$ to each vector in $T_mM$, the usual degree to each homogeneous operator on $T^*_mM$,   degree $0$ to each $dx^i$ and degree $1$ to each $d \xi_j$.  For instance, the terms in  \eqref{eq-example-of-symbol-operator} have degrees $1$ and $0$, respectively.

\begin{lemma} 
\label{lem-symbol-homomorphism-on-pdo}
Let $\nabla$ be a torsion-free connection on $M$
and let  $m\in M$. There is a unique isomorphism of complex $\Z$-graded algebras 
\[
\Symbol^\nabla _m \colon  \gr \PDO (M,S)  \big \vert _m
  \stackrel \cong  \longrightarrow  \SymbolAlg_m
\]
such that for all local coordinates near $m$, 
\begin{enumerate}[\rm (i)]

\item    $[\nabla_{(\partial/\partial x^i)^H}]\mapsto \partial/\partial x^i  $ and  $[\nabla _{\partial/\partial \xi _j}]\mapsto \partial /\partial \xi _j$.

\item If $f$ is any polyhomogeneous function on $\Tcirc M$, then $f\mapsto  f\vert _{\Tcirc _m M} $.

\item       $ [\psi^i] \mapsto  dx^i$ and $[\overline \psi _j]\mapsto d\xi_j$.

\end{enumerate}
\end{lemma} 

\begin{proof} 
The reversed correspondences define a morphism of graded algebras in the reverse direction by the universal property inherent in the definition of $\SymbolAlg_m$  as a tensor product.   It follows from the  local coordinate description of the Perrot order in Lemma~\ref{lem-order-bound-by-coordinate-filtration} that this morphism is an isomorphism.
\end{proof} 

\begin{remark}
Even though the Perrot order and  therefore  the associated graded algebra $\gr \PDO (M,S)$, are independent of the choice of   connection used in their definitions, the  homomorphism in Lemma~\ref{lem-symbol-homomorphism-on-pdo} \emph{does} depend on the choice of connection.  
\end{remark}

We shall now extend the symbol morphism $\Symbol^\nabla_m$ from polyhomogeneous differential operators to elements of the bimodule $\FB(M,S)$.  To do so we shall first look to Definition~\ref{def-definition-of-fbu} to find a suitable target for the symbol homomorphism.

\begin{definition} 
Define the \emph{bimodule order} of an element in $A_m$ by 
\[
\begin{aligned}
\bimodorder \Bigl (\frac{\partial}{\partial x^i}\Bigr ) = 3,\quad 
& \bimodorder \Bigl (\frac{\partial}{\partial \xi_j}\Bigr ) = -1 \\
\bimodorder \bigl( \psi^i \bigr) = 0 ,\quad 
& \bimodorder \bigl (\overline \psi_j\bigr ) = 0
\end{aligned}
\]
and $\bimodorder (f) = 2 \degree (f)$ for a homogeneous function $f$.
Denote by $\FA_m$  the space of formal Laurent series 
\[
\sum \varepsilon ^k D_k \exp(\varepsilon \Delta)  \in A_m [\varepsilon ^{-1} , \varepsilon ]]
\]
with finite singular parts for each of  which there exists $N$ with 
\[
\bimodorder (D_k) \le k + N 
\]
 for all $k$.   
\end{definition}

The space $\FA_m$ is a bimodule over $\SymbolAlg_m$, and applying the symbol homomorphism termwise to a Laurent series, we obtain a morphism
\[
\Symbol_m^\nabla \colon 
\gr \FB (M,S) \vert _m  \longrightarrow  \FA_m
\]
that is compatible with bimodule structures. 
We are  going to factor Perrot's supertrace on order zero elements through this symbol morphism, as follows: we shall define a supertrace morphism 
\[
\str_m \colon \FA _m \longrightarrow \Density (M)\vert _m 
\]
and show that 
\[
\str(X)\vert _m = \str_m (\Symbol^\nabla_m (X))
\]
for every $X \in \FB (M,S)$ of Perrot order $0$. 

To begin the definition of $\str_m$, suppose we are given an order $-n$ homogeneous smooth function  $f$ on $\Tcirc_m M$, and let $\alpha \in \Lambda^n T_mM$.  Think of $\alpha$ as a (translation-invariant) top-degree form on $\Tcirc _m M$.  If $E$ is the Euler vector field, as usual, then the contracted  form $\iota _E (f \cdot \alpha)$ is basic for the action of the positive real numbers on $\Tcirc _mM$ by scalar multiplication, and it is therefore the pullback of a form $\sigma_{f, \alpha}$ on the sphere $S^*_m M$.  Consider next the integral
 \begin{equation}
 \label{eq-integral-of-sigma-f-alpha}
 \int _{S^*_mM} \sigma_{f , \alpha} ,
 \end{equation}
 formed using the orientation on the sphere for which the integral is positive when $f$ is positive. This orientation depends on $\alpha$, and so the integral is not bilinear in $f$ and $\alpha$, but rather 
  \[
 \int _{S^*_mM} \sigma_{f, t \alpha} = |t|  \int _{S^*_mM} \sigma_{f, \alpha}\qquad \forall t \in \R.
 \]
 The map $\alpha \mapsto \int _{S^*_mM} \sigma_{f,  \alpha}$ therefore defines a density at $m$ (see \cite[pp.30-31]{BGV} or \cite[Ch.14]{Lee} again), which we shall write as 
 \begin{equation}
     \label{eq-def-of-fint-at-m}
 \fint _{S^* _m M} f\in \Density (M)\vert _m 
 \end{equation}
(and as before we extend this to polyhomogeneous $f$ by selecting the degree $-n$ component). This construction has the property that if $f$ is a smooth, degree $-n$ function on $\Tcirc M$ that is compactly supported in the $M$-direction, then 
 \begin{equation}
     \label{eq-fubini-density-formula}
 \int _M \Bigl (  \fint_{S^*_mM} f\vert_{\Tcirc _mM}\Bigr ) = \fint _{S^*M} f .
 \end{equation}
 
 To define a trace\footnote{The functional is indeed a trace, although we shall not prove this fact because we shall not need it.}  functional on $\FA_m$, we  use the integral above and essentially the same contraction operation that we used for Perrot's trace: namely for $D \in \SymbolAlg_m$ we define 
 \[
 \bigl \langle D \exp (\varepsilon \Delta )\bigr \rangle = 
\varepsilon^{-n} \bigl ( D \exp (\varepsilon ^{-1} (x^i-y^i)(\xi_i-\eta_i)) \bigr ) \Big \vert _{x^i=y^i=0,\,\, \xi_i=\eta_i} .
 \]
 This is a polyhomogeneous function on $\Tcirc _m M$ with values in the exterior algebra $\Lambda^*(T_mM \oplus T^*_m M)$.  Then we define 
 \begin{multline}
 \label{eq-def-of-str-m}
 \str_{m,F} \bigl ( D \exp (\varepsilon \Delta )\bigr ) =
 \text{coefficient of}
 \\
 dx^1\cdots dx^nd\xi_1\cdots d \xi _n \,\, \text{in} \,\,
 \fint_{S^*_mM} \bigl \langle D \exp (\varepsilon \Delta ) \bigr \rangle ,
 \end{multline}
 and then extend termwise to the Laurent series in $\FA_m$.  The result is a Laurent series in $\varepsilon$ and, just as we did when we constructed Perrot's trace, we then define $\str_m$ by taking the coefficient of $\varepsilon^0$.

Thanks to our direct mimicry of the constructions in Section~\ref{sec-perrot-trace}, and  thanks to \eqref{eq-fubini-density-formula}, the following result is clear: 

\begin{lemma} 
Let $x_1,\dots, x_n$ be local coordinates on $M$, defined near a point $m\in M$. If $\nabla$ is the canonical flat connection on $TM$ associated to these coordinates \textup{(}in which all the Christoffel symbols $\Gamma^k_{ij}$ are identically zero\textup{)},  then 
\[
 \str (X)\big \vert_m   = \str_m (\Symbol_m^\nabla (X))
\]
for every $X \in \FB(M,S)$ of Perrot order $0$. \qed
\end{lemma} 

In fact, the same result holds for \emph{any} connection: 
 
\begin{theorem} 
\label{thm-general-connection-formula-for-trace-density}
If $\nabla$ is \emph{any} torsion-free connection on $TM$, if $m\in M$, and if 
\[
\Symbol_m^\nabla \colon \FB_0(M)  \longrightarrow \FA _m
\]
is the symbol map at $m$ defined by $\nabla$ and acting on Perrot order zero operators,  then 
\[
 \str (X)\big \vert_m   = \str_m (\Symbol_m^\nabla (X))
\]
for every $X \in \FB^0(M)$.
\end{theorem}  

\begin{proof} 
Given $\nabla$ and $m$, choose local coordinates near $m$ for which all the Christoffel symbols for $\nabla$ vanish at $m$.  Then define $\nabla^0$ to be the canonical flat connection for these coordinates.  By examining the two morphisms $\Symbol_m^\nabla$ and $\Symbol_m^{\nabla^0}$  in this given coordinate system, we see that these two symbol morphisms are in fact equal.  So the theorem follows from the previous lemma.
\end{proof} 

\begin{remark} Needless to say, the theorem relies heavily on the coordinate independence of Perrot's trace (and in particular the theorem is far from trivial).
\end{remark}

\section{Computation of the JLO cocycle} 
\label{sec-computation-of-jlo-cocycle}

In this final section we shall use Theorem~\ref{thm-general-connection-formula-for-trace-density} and one final additional computation to identify  the cocycle  $JLO^{\Dirac, \delta}$ from Section~\ref{sec-jlo-type-cocycle} with the cocycle associated with the Todd class that was described in Section~\ref{sec-introduction}.  The starting point is the following computation:

\begin{lemma} If $\nabla$ is any torsion-free connection on $M$, then 
\[
\perrotorder\bigl (\Dirac_{\mathrm{vert}}^2\bigr )  \le 0
\quad
\text{and} 
\quad
\perrotorder\bigl ( \{ \Dirac_{\mathrm{horiz}}, \Dirac_{\mathrm{vert}} \}\bigr ) \le 0 .
\]
In addition if $a$ is any smooth, degree zero function on $\Tcirc M$, then 
\[
\perrotorder\bigl ([\Dirac  ,a] \bigr ) \le 0 .
\]
\end{lemma}

\begin{proof} 
This evident from the formulas in Definition~\ref{def-d-vert-and-d-horiz} and in Lemmas~\ref{lem-square-of-d-horiz} and \ref{lem-local-formula-for-cross-term}.
\end{proof} 

It follows that the argument of the supertrace in the formula for $JLO^{\Dirac,\delta}(a^0,\dots, a^p)$ has Perrot order $0$, and therefore we may use the  techniques of Section~\ref{sec-perrot-order} to compute the supertrace.

With this in mind, we  now fix a point $m\in M$.  As explained in the Section~\ref{sec-perrot-order}, the quantity $JLO^{\Dirac,\delta}(a^0,\dots, a^p)$ is the integral of a density on $M$; we shall denote by 
\[
JLO^{\Dirac,\delta}(a^0,\dots, a^p)\vert_m \in \Density(M)\vert _m
\]
the value of this density at $m$. 

If we  write 
\[
R \Bigl(\frac{\partial}{\partial x^i}, \frac{\partial}{\partial x^j}\Bigr )\, \frac{\partial}{\partial x^k} =   R^k_{ij\ell}\, \frac {\partial}{\partial x^\ell} ,
\]
as we did earlier, and then define 
\[
R_{\ell}^k = \tfrac 12 R^k_{ij\ell}dx^idx^j \big \vert_m \in \Lambda ^2 T_m^*M ,
\]
then 
\begin{equation}
    \label{eq-sybol-formula1}
\Symbol^\nabla_m (\exp (s \Dirac^2))  = \exp \bigl ( s (\varepsilon  \Delta {+} \varepsilon^2 \xi_k R^k_\ell {\partial}/{\partial \xi _\ell}))\bigr )
\end{equation}
Now  choose coordinates $x^1,\dots, x^n$ near $m$ for  which the Christoffel symbols $\Gamma_{ij}^k$ vanish at $m$. Then 
\begin{equation}
    \label{eq-sybol-formula2}
\begin{gathered}
\Symbol^\nabla_m ([\Dirac, a])   = \varepsilon \frac{\partial a}{\partial x^i} dx^i + \frac{\partial a}{\partial \xi_j} d\xi _j 
\\
\Symbol^\nabla_m (\delta(\Dirac))   =  -   \varepsilon q^{-1}\frac{\partial q}{\partial x^i} dx^i - q^{-1} \frac{\partial q}{\partial \xi_j} d\xi _j 
\\
\end{gathered}
\end{equation}
where as in Section~\ref{sec-jlo-type-cocycle}, $q\colon \Tcirc M\to (0,\infty)$ is the norm-function associated to a metric on $T^*M$. 

For brevity, we shall write 
\[
\begin{gathered}
\Symbol^\nabla_m ([\Dirac, a]) = d_\varepsilon a = \varepsilon d_{\mathrm{horiz}}a + d_{\mathrm{vert}}a  
\\
\Symbol^\nabla_m (\delta(\Dirac))  = - q^{-1} d_\varepsilon q =   -   \varepsilon q^{-1}d _{\mathrm{horiz}} q - q^{-1} d_{\mathrm{vert}}q .
\\
\end{gathered}
\]
In addition, to streamline some of the formulas that follow, we shall write 
\[
\Delta _R  =  \Delta {+} \varepsilon \xi_k R^k_{\ell} {\partial}/{\partial \xi _\ell} =  \Delta {+} \varepsilon \xi \cdot  R \cdot {\partial}/{\partial \xi} .
\]
Then according to the formulas \eqref{eq-sybol-formula1} and \eqref{eq-sybol-formula2} above and the formula in Theorem~\ref{thm-jlo-type-cocycle}, the density  $JLO^{\Dirac,\delta}(a^0,\dots, a^p)\vert_m  $   is the sum from $k=1$ to $k={p{+}1}$ of the following integrals:
\begin{multline}
\label{eq-pointwise-form-of-jlo-d-delta-1}
(-1)^{k-1}\int_{\Sigma^{p+1}}
 \str_m \bigl ( a^0 \exp   ( s_0 \varepsilon \Delta_R)  d_\varepsilon a^1 \exp   ( s_1 \varepsilon \Delta_R) \cdots 
 \\
 \dots d_\varepsilon a^{k-1} \exp   ( s_{k-1} \varepsilon \Delta_R) q^{-1} dq \exp   ( s_{k} \varepsilon \Delta_R) d_\varepsilon a^{k} \cdots \\
 \cdots \exp   ( s_{p-1} \varepsilon \Delta_R) d_\varepsilon a^{p-1}   \exp   ( s_p \varepsilon \Delta_R)\bigr ) \, ds  .
\end{multline}
Using the formula 
\[\exp (s \varepsilon \Delta_R ) X \exp (-s \varepsilon \Delta_R ) = \exp ( s \varepsilon \ad_{\Delta_R})(X)
\]
we can write the argument of $\str_m$ in \eqref{eq-pointwise-form-of-jlo-d-delta-1} as an infinite linear combination (convergent in the adic topology) of terms
\begin{equation}
\label{eq-pointwise-form-of-jol-d-delta-expanded}
\varepsilon ^k X_0 \ad_{\Delta_R}^{k_0}(X_1)
 \ad_{\Delta_R}^{k_1}(X_2)
 \cdots 
  \ad_{\Delta_R}^{k_{p-2}}(X_{p-1}) \ad_{\Delta_R}^{k_{p-1}}(X_p) \exp (\varepsilon \Delta _R) ,
\end{equation}
 where each $X_j$ is one of the $d_\varepsilon a^i$ or $q^{-1} d_\varepsilon q$.  W

 Next,  as long as any $X\in \SymbolAlg_m$ includes no $\xi$-derivatives
\[
\ad_{\Delta_R}(X) = \bigl (\partial/\partial x^j + \varepsilon \xi_i R_{ij}\bigr ) [\partial/\partial \xi_j , X],
\]
and we note that in this case  the result on the right-hand side also includes no $\xi$-derivatives.  Now we appeal to the following result, whose proof we shall defer for a moment: 

\begin{proposition}
\label{prop-schur-calculation1}
If an element 
$
X\in A_m
$
 includes no $\xi$-derivatives, then 
\[
\str_m \bigl (\varepsilon ^ r (\partial_{x^j} + \varepsilon \xi_i R_{ij})  X \exp (\varepsilon \Delta_R   )  \bigr )
 = 0 
 \]
for all $r$ and all $j$. 
\end{proposition}

It follows from the lemma that all the terms \eqref{eq-pointwise-form-of-jlo-d-delta-1} have vanishing trace, except when all $k_j$ are zero, and hence that \eqref{eq-pointwise-form-of-jlo-d-delta-1} is equal to
\begin{equation*}
(-1)^{k-1} \int_{\Sigma^{p+1}}
 \str_m \bigl (  a^0    d_\varepsilon a^1   \cdots d_\varepsilon a^{k-1} q^{-1} d_\varepsilon q \cdots  d_\varepsilon a^ p    \exp   ( \varepsilon  \Delta_R)\bigr )   \, ds .
\end{equation*}
The integrand now no longer depends on $s$, and so the integral is simply 
\[
 \frac{(-1)^{k-1}}{(p{+}1)!} \str_m \bigl (   a^0    d_\varepsilon a^1   \cdots d_\varepsilon a^{k-1} q^{-1} d_\varepsilon q \cdots  d_\varepsilon a^ p    \exp   ( \varepsilon  \Delta_R) \bigr ).
\]
In summary, then: 
\begin{multline}
\label{eq-pointwise-form-of-jlo-d-delta-2}
JLO^{\Dirac,\delta}(a^0,\dots, a^p)\vert_m 
  \\
    \begin{aligned}
  & =  \sum_{k=1}^{p+1}  \frac{(-1)^{k-1}}{(p{+}1)!}   
  \str_m \bigl (   a^0    d_\varepsilon a^1 \cdots  d_\varepsilon a^{k-1} q^{-1} d_\varepsilon q\cdots  d_\varepsilon a^ p    \exp   ( \varepsilon \Delta_R)  \bigr ) 
 \\
 &\qquad  = \sum_{k=1}^{p+1} \frac{1}{(p{+}1)!}  \str_m \bigl (  q^{-1} d_\varepsilon q\,
   a^0    d_\varepsilon a^1   \cdots  d_\varepsilon a^ p    \exp   (  \varepsilon \Delta_R) \bigr )   
 \\
 & \qquad \qquad = \frac{1}{p!}\str_m \bigl ( q^{-1} d_\varepsilon q\,   a^0    d_\varepsilon a^1   \cdots  d_\varepsilon a^ p    \exp   (  \varepsilon \Delta_R)\bigr ) .
 \end{aligned}
\end{multline}
Now we quote another result, whose proof, like that of Proposition~\ref{prop-schur-calculation1}, we shall   defer for a short time:

 \begin{proposition} 
  \label{prop-schur-calculation2}
If  $X \in A_m$ has no $x$- or $\xi$-derivatives, then    
\[
  \bigl \langle   X \exp (\varepsilon \Delta_R   )  \bigr \rangle 
 =  \varepsilon ^{-n} X \Todd (\varepsilon ^2 R)  .
 \]
\end{proposition}

To continue the computation, let us now write  
\begin{multline*}
\left [ q^{-1} d_\varepsilon q\,   a^0    d_\varepsilon a^1   \cdots  d_\varepsilon a^ p \Todd (\varepsilon ^2 R)\right ]^{\mathrm{top}}
\\
=\text{top exterior form-degree part of}
\\
q^{-1} d_\varepsilon q\,   a^0    d_\varepsilon a^1   \cdots  d_\varepsilon a^ p \Todd (\varepsilon ^2 R).
\end{multline*}
The Todd class is a polynomial in $\varepsilon$, for which the coefficient of $\varepsilon^k$ lies in $\Lambda ^k T^*M$ (with only even $k$ appearing). Each of the $d_\varepsilon$-differentials lies in $\varepsilon \Lambda^1 T^*_mM + \Lambda^1 T_mM$.   It therefore follows from an examination of  degree in $\Lambda ^* T^*_mM$  that 
\begin{multline*}
\left [ q^{-1} d_\varepsilon q\,  a^0  d_\varepsilon a^1    \cdots  d_\varepsilon a^ p \Todd (\varepsilon ^2 R)\right ]^{\mathrm{top}} 
\\
=
\varepsilon ^n \left [ q^{-1} d  q\,   a^0   da^1     \cdots  d  a^ p \Todd (  R)\right ]^{\mathrm{top}}.
\end{multline*}
As a result, it follows from \eqref{eq-pointwise-form-of-jlo-d-delta-2} and  Proposition~\ref{prop-schur-calculation2}, and the definition of $\str_m$   in \eqref{eq-def-of-str-m}  that
\begin{multline}
\label{eq-pointwise-form-of-jlo-d-delta-3}
   \frac{1}{p!}\str_m \bigl ( q^{-1} d_\varepsilon q\,   a^0    d_\varepsilon a^1   \cdots  d_\varepsilon a^ p    \exp   (  \varepsilon \Delta_R)\bigr )
     \\
    = \text{coefficient of}\,\, 
 dx^1\cdots dx^nd\xi_1\cdots d \xi _n 
 \\
 \text{in} \,\,\, 
 \frac{1}{p!} \fint_{S^*_m M} [ q^{-1} d  q\,  a^0 da^1    \cdots  d  a^ p \Todd ( R) ]^{\mathrm{top}} .
\end{multline}
In this formula $[ q^{-1} d  q\,  a^0 da^1    \cdots  d  a^ p \Todd ( R) ]^{\mathrm{top}}$   is to be regarded as a polyhomogeneous, exterior algebra-valued function on $\Tcirc_mM$ and then written as a scalar function times $dx^1\cdots dx^nd\xi_1\cdots d \xi _n$;  the integral $\fint$ is applied to this scalar function as  in \eqref{eq-def-of-fint-at-m}.  

The   right-hand side in \eqref{eq-pointwise-form-of-jlo-d-delta-3} involves a contraction by the Euler vector field, so let us now note that 
\[
\begin{aligned}
\iota_E \bigl ( q^{-1} d  q\,  a^0 da^1    \cdots  d  a^ p \Todd ( R) \bigr ) 
& = 
\iota_E \bigl ( q^{-1} d  q\bigr )    a^0  da^1   \cdots  d  a^ p \Todd ( R)   \\
& =   a^0  da^1   \cdots  d  a^ p \Todd ( R) ,
\end{aligned}
\]
since $\iota_E \bigl (  q^{-1} d  q\bigr ) =1$ while $\iota_E\bigl ( a^0  da^1   \cdots  d  a^ p \Todd ( R)\bigr ) =0$.   
It follows that
\begin{multline*}
\str_m \bigl ( q^{-1} d_\varepsilon q\,   a^0    d_\varepsilon a^1   \cdots  d_\varepsilon a^ p    \exp   (  \varepsilon \Delta_R)\bigr )
\\
= \pi_*\bigl (  [a^0 da^1\cdots da^p\Todd(R)]^{\mathrm{top}}\bigr ) \big \vert_m ,
\end{multline*}
where $\pi \colon S^*M \to M$ and where the differential form is treated as a density using the orientation of $S^*M$ given in Definition~\ref{def-orientation-on-s-star-m}.
We have arrived at our goal: 

 \begin{theorem} 
Let $M$ be a smooth, closed manifold, let $\Delta $ be a torsion-free connection on $TM$ with curvature $R$, and let $\Dirac$ be the associated Dirac operator, as in Definition~\textup{\ref{def-dirac-operator}}.  If $S^*M$ is equipped with the orientation in Definition~\textup{\ref{def-orientation-on-s-star-m}}, then
\[
JLO^{\Dirac, \delta}(a^0, \dots, a^p) 
= \frac{1}{p!}  
\int _{S^*M} a^0  da^1    \cdots da^ p
  \Todd  (  R   )   
\]
for all $p>0$ and all  $a^0,\dots, a^p\in C^\infty (S^*M)$. \qed
 \end{theorem} 
 
 It remains to prove Propositions~\ref{prop-schur-calculation1} and \ref{prop-schur-calculation2}.

 \begin{proof}[Proof of Propositions~\ref{prop-schur-calculation1} and \ref{prop-schur-calculation2}]
 The following beautiful argument (along with everything else in this section) is due to Perrot \cite[Lem.4.4]{Perrot2013}, to which we refer for full details.  It is reminiscent of approaches to the Baker-Campbell-Haus\-dorff formula going back to Schur
 \cite{Schur1889}.
 
 For $t\in \R$ form the element 
 \[
 Y(t) = \exp (\varepsilon \Delta + t \varepsilon^2 \xi \cdot  R \cdot {\partial}/{\partial \xi} )\exp (-\varepsilon \Delta)
 \]
 in $\SymbolAlg_m[\varepsilon]$; it is a polynomial function of $t$ and $\varepsilon $ thanks to the nilpotence of the subspace  
 \[
 \Lambda ^2(T_mM\oplus T^*_mM)\subseteq \Lambda ^*(T_mM\oplus T^*_mM).
 \]
 If $X\in \SymbolAlg_m[\varepsilon]$, then according to \eqref{eq-reconciliation-with-perrot-contraction},
 \begin{multline*}
 \bigl \langle X \exp ( \varepsilon \Delta + t \varepsilon^2 \xi \cdot  R \cdot {\partial}/{\partial \xi} )\bigr \rangle 
 \\
 =
 \varepsilon^{-n} \Bigl ( X\cdot  Y(t) \exp \bigl (-\varepsilon ^{-1} (x^i-y^i)(\xi_i-\eta_i)\bigr ) \Bigr ) \Big \vert _{x^i=y^i,\,\, \xi_i=\eta_i}
 \end{multline*}
 Consider now the (form-valued) function 
 \[
 H(t) = \varepsilon^{-n} 
 Y(t) \exp \bigl (-\varepsilon ^{-1} (x^i-y^i)(\xi_i-\eta_i)\bigr )
 \]
 on the product of $\Tcirc T_mM$ with itself.  A direct computation shows that it satisfies a differential equation of the type 
 \[
 \frac{dH}{dt}(t)  = L(t)  H(t) ,
 \]
 where $L(t)$ is a polynomial in $t$ and $\varepsilon$ with coefficients in the algebra $\SymbolAlg_m$, and that the function 
 \[
 \Todd (\varepsilon ^2 t R) \cdot \exp \Bigl ( -\eta \cdot \varepsilon  R \cdot (x{-}y) -   (\xi {-} \eta) \cdot 
 \frac{t \varepsilon R}{1 -\exp (- t \varepsilon^2R )}\cdot (x {-}y)
 \Bigr ) 
 \]
 satisfies the same differential equation, with the same initial condition at $t=0$.  The two solutions are necessarily equal, and the two propositions follow by setting $t=1$ in the explicit solution, applying respectively $  (\partial/\partial x^j + \varepsilon \xi_i R_{ij}  ) X$, as in the first proposition, or  $X$ alone   
 as in the second, and restricting to the diagonal to obtain an explicit formula for the contractions.
 \end{proof}

\bibliography{References}

\begin{thebibliography}{H{o}r83}

\bibitem[BGV92]{BGV}
Nicole Berline, Ezra Getzler, and Mich\`ele Vergne.
\newblock {\em Heat kernels and {D}irac operators}, volume 298 of {\em
  Grundlehren der mathematischen Wissenschaften}.
\newblock Springer-Verlag, Berlin, 1992.

\bibitem[CM95]{ConnesMoscovici95}
Alain Connes and Henri Moscovici.
\newblock The local index formula in noncommutative geometry.
\newblock {\em Geom. Funct. Anal.}, 5(2):174--243, 1995.

\bibitem[CM98]{ConnesMoscovici98}
Alain Connes and Henri Moscovici.
\newblock Hopf algebras, cyclic cohomology and the transverse index theorem.
\newblock {\em Comm. Math. Phys.}, 198(1):199--246, 1998.

\bibitem[Con94]{ConnesNCGBook}
Alain Connes.
\newblock {\em Noncommutative geometry}.
\newblock Academic Press, Inc., San Diego, CA, 1994.

\bibitem[Get83]{Getzler83}
Ezra Getzler.
\newblock Pseudodifferential operators on supermanifolds and the
  {A}tiyah-{S}inger index theorem.
\newblock {\em Comm. Math. Phys.}, 92(2):163--178, 1983.

\bibitem[Hig04]{Higson04}
Nigel Higson.
\newblock The local index formula in noncommutative geometry.
\newblock In {\em Contemporary developments in algebraic {$K$}-theory}, ICTP
  Lect. Notes, XV, pages 443--536. Abdus Salam Int. Cent. Theoret. Phys.,
  Trieste, 2004.

\bibitem[H{o}r83]{HormanderVol1}
Lars H{o}rmander.
\newblock {\em The analysis of linear partial differential operators. {I}},
  volume 256 of {\em Grundlehren der mathematischen Wissenschaften}.
\newblock Springer-Verlag, Berlin, 1983.
\newblock Distribution theory and Fourier analysis.

\bibitem[HY19]{HigsonYi19}
Nigel Higson and Zelin Yi.
\newblock Spinors and the tangent groupoid.
\newblock {\em Doc. Math.}, 24:1677--1720, 2019.

\bibitem[JLO88]{JLO88}
Arthur Jaffe, Andrzej Lesniewski, and Konrad Osterwalder.
\newblock Quantum {$K$}-theory. {I}. {T}he {C}hern character.
\newblock {\em Comm. Math. Phys.}, 118(1):1--14, 1988.

\bibitem[Kas89]{Kassel89}
Christian Kassel.
\newblock Le r\'{e}sidu non commutatif (d'apr\`es {M}. {W}odzicki).
\newblock {\em Ast\'{e}risque}, (177-178):Exp. No. 708, 199--229, 1989.
\newblock S\'{e}minaire Bourbaki, Vol. 1988/89.

\bibitem[Lee03]{Lee}
John~M. Lee.
\newblock {\em Introduction to smooth manifolds}, volume 218 of {\em Graduate
  Texts in Mathematics}.
\newblock Springer-Verlag, New York, 2003.

\bibitem[Nis97]{Nistor97}
Victor Nistor.
\newblock Higher index theorems and the boundary map in cyclic cohomology.
\newblock {\em Doc. Math.}, 2:263--295, 1997.

\bibitem[Per12]{Perrot12}
Denis Perrot.
\newblock Extensions and renormalized traces.
\newblock {\em Math. Ann.}, 352(4):911--940, 2012.

\bibitem[Per13a]{Perrot13b}
Denis Perrot.
\newblock Local index theory for certain {F}ourier integral operators on {L}ie
  groupoids.
\newblock Preprint, 2013.
\newblock \href{http://arxiv.org/abs/1401.0225}{arXiv:1401.0225}.

\bibitem[Per13b]{Perrot2013}
Denis Perrot.
\newblock Pseudodifferential extension and {T}odd class.
\newblock {\em Adv. Math.}, 246:265--302, 2013.

\bibitem[Per16]{Perrot16}
Denis Perrot.
\newblock Index theory for improper actions: localization at units.
\newblock Preprint, 2016.
\newblock \href{http://arxiv.org/abs/1612.04090}{arXiv:1612.04090}.

\bibitem[PR14]{PerrotRodsphon14}
Denis Perrot and Rudy Rodsphon.
\newblock An equivariant index theorem for hypoelliptic operators.
\newblock Preprint, 2014.
\newblock \href{http://arxiv.org/abs/1412.5042}{arXiv:1412.5042}.

\bibitem[Qui88]{Quillen88}
Daniel Quillen.
\newblock Algebra cochains and cyclic cohomology.
\newblock {\em Inst. Hautes \'{E}tudes Sci. Publ. Math.}, (68):139--174 (1989),
  1988.

\bibitem[Rad91]{Radul91}
Andrey~O. Radul.
\newblock Lie algebras of differential operators, their central extensions and
  {$W$}-algebras.
\newblock {\em Funktsional. Anal. i Prilozhen.}, 25(1):33--49, 1991.

\bibitem[Rod15]{Rodsphon2015}
Rudy Rodsphon.
\newblock Zeta functions, excision in cyclic cohomology and index problems.
\newblock {\em J. Funct. Anal.}, 268(5):1167--1204, 2015.

\bibitem[Roe98]{RoeAsymptoticMethods}
John Roe.
\newblock {\em Elliptic operators, topology and asymptotic methods}, volume 395
  of {\em Pitman Research Notes in Mathematics Series}.
\newblock Longman, Harlow, second edition, 1998.

\bibitem[Sch89]{Schur1889}
Friedrich Schur.
\newblock Neue {B}egr\"{u}ndung der {T}heorie der endlichen
  {T}ransformationsgruppen.
\newblock {\em Math. Ann.}, 35(1-2):161--197, 1889.

\bibitem[Wod87]{Wodzkicki87}
Mariusz Wodzicki.
\newblock Noncommutative residue. {I}. {F}undamentals.
\newblock In {\em {$K$}-theory, arithmetic and geometry ({M}oscow,
  1984--1986)}, volume 1289 of {\em Lecture Notes in Math.}, pages 320--399.
  Springer, Berlin, 1987.

\end{thebibliography}
\bibliographystyle{alpha}

\end{document}